\documentclass[11pt,a4wide]{article}

\usepackage{amsmath}
\usepackage{amsthm}
\usepackage{amsxtra}
\usepackage{amsfonts,amssymb}
\usepackage[dvips]{graphics}
\usepackage{eepic}
\newtheorem{Th}{Theorem}
\newtheorem{Prop}[Th]{Proposition}
\newtheorem{Lm}[Th]{Lemma}
\newtheorem{Co}[Th]{Corollary}

\theoremstyle{definition}
\newtheorem{Def}[Th]{Definition}

\newtheorem{Rem}{Remark}

\renewcommand{\theparagraph}{\thesection}

\numberwithin{equation}{section}

\date{}

\begin{document}
\title{
\begin{center}
A description of characters on the infinite wreath product.
\end{center}}

\author{A.V. Dudko, N.I. Nessonov\footnote{Supported by
the CRDF-grant UM1-2546}}

\maketitle

\begin{abstract}
Let $\mathfrak{S}_\infty$ be the infinity permutation group and $\Gamma$ an
arbitrary group. Then $\mathfrak{S}_\infty$ admits a natural action on
$\Gamma^\infty$ by automorphisms, so one can form a semidirect product
$\Gamma^\infty\rtimes \mathfrak{S}_\infty$, known as the {\it wreath}
product $\Gamma\wr\mathfrak{S}_\infty$ of $\Gamma$ by
$\mathfrak{S}_{\infty}$. We obtain a full description of unitary
$I\!I_1-$factor-representations of $\Gamma\wr\mathfrak{S}_\infty$ in terms
of finite characters of $\Gamma$. Our approach is based on extending
Okounkov's classification method for admissible representations of
$\mathfrak{S}_\infty\times\mathfrak{S}_\infty$. Also, we discuss certain
examples of representations of type $I\!I\!I$, where the {\it modular
operator} of Tomita-Takesaki expresses naturally by the asymptotic
operators, which are important in the characters-theory of
$\mathfrak{S}_\infty$.

\end{abstract}

\section{Introduction}

\paragraph{\theparagraph.1. A definition of the wreath product.}
Let $\mathbb{N}$ stand for the natural numbers. A bijection
$s:\mathbb{N}\to\mathbb{N}$ is called {\it finite} if the set
$\left\{i\in\mathbb{N}|s(i)\ne i\right\}$ is finite. Define
$\mathfrak{S}_\infty$ as the group of all finite bijections
$\mathbb{N}\to\mathbb{N}$ and set
$\mathfrak{S}_n=\left\{s\in\mathfrak{S}_\infty|\;s(i)=i\;\forall\;i>n
\right\}$. For every group $\Gamma$, an element of $\Gamma^n$ can always be
written as a sequenced collection $\left[\gamma_k\right]_{k=1}^n=
\left(\gamma_1,\gamma_2,\ldots,\gamma_n\right)$, where $\gamma_k\in\Gamma$.
Let $e$ be the unit of $\Gamma$. For any $n>1$ we identify the element
$\left(\gamma_1,\gamma_2,\ldots,\gamma_{n-1}\right)\in\Gamma^{n-1}$ with
$\left(\gamma_1,\gamma_2,\ldots,\gamma_{n-1},e\right)\in\Gamma^n$ and set
$\Gamma^\infty_e=\underrightarrow{\lim}\Gamma^n$. One can view
$\Gamma^\infty_e$ as a group of infinite sequenced collections
$\left[\gamma_k\right]_{k=1}^\infty$ such that there are finitely many
elements $\gamma_k$ not equal to $e$. The {\it wreath} product
$\Gamma\wr\mathfrak{S}_n$ is the semidirect product
$\Gamma^n\rtimes\mathfrak{S}_n$ for the natural permutation action of
$\mathfrak{S}_n$ on $\Gamma^n$ (see \cite{EM}). In the same way, we define
the group $\Gamma\wr\mathfrak{S}_\infty=\Gamma^\infty_e
\rtimes\mathfrak{S}_\infty$. $\Gamma\wr\mathfrak{S}_\infty$ can be also
viewed as the inductive limit
$\underrightarrow{\lim}\Gamma\wr\mathfrak{S}_n$. Using the embeddings
$\gamma\in\Gamma^n\to\left(\gamma,id\right)\in\Gamma\wr\mathfrak{S}_n$ and
$s\in\mathfrak{S}_n\to\left(e^{(n)},s\right)\in\Gamma\wr\mathfrak{S}_n$,
where $e^{(n)}=\left(e,e,\ldots,e\right)$ and $id$ is the identical
bijection, we may identify $\Gamma^n$ and $\mathfrak{S}_n$ with the
corresponding subgroups in $\Gamma\wr\mathfrak{S}_n$. If $\Gamma$ is a
topological group, then we equip $\Gamma^n$ with the natural
product-topology. Furthermore, we will always consider $ \Gamma^\infty_e$
as a topological group with the inductive limit topology. As a set,
$\Gamma\wr\mathfrak{S}_\infty$ is just
$\Gamma^\infty_e\times\mathfrak{S}_\infty$. Therefore, we equip
$\Gamma\wr\mathfrak{S}_\infty$ with the product-topology, considering
$\mathfrak{S}_\infty$ as a discrete topological space.

\medskip
\paragraph{\theparagraph.2. The Results.}
In this paper we give a full classification of {\it indecomposable}
characters (see Definitions (\ref{trace}) -- (\ref{def4})) on
$\Gamma\wr\mathfrak{S}_\infty$ (Theorem {\ref{mainth}}). Our approach is
based on the semigroup method of Olshanski \cite{O2} and the ideas of
Okounkov used in the study of {\it admissible} representations of the
groups related to $\mathfrak{S}_\infty$\cite{Ok1},\cite{Ok2}. We have
noticed that two double cosets containing the transposition or
$\gamma\in\Gamma$ are commutated, as the elements of Olshanski semigroup.
(see Fig. \ref{pic8}, p.\pageref{com} and Lemma \ref{abelian}). This
observation enables one to develop Okounkov's method for the group
$\mathfrak{S}_\infty\wr\Gamma$ (see Section {\ref{mainresult}}). In Section
{\ref{examples}} we discuss certain examples of representations of type
$I\!I\!I$. The corresponding positive definite functions (p.d.f.) $\varphi$
are not characters, but the following holds:
\begin{eqnarray}\label{kms}
\varphi(sg)=\varphi(gs)\textit{ for all
}g\in\mathfrak{S}_\infty\wr\Gamma\textit{ and }s\in\mathfrak{S}_\infty.
\end{eqnarray}
Hence the restriction $\varphi\Big|_{\mathfrak{S}_\infty}$ is a character.
At that, the Okounkov's asymptotic operators (see (\ref{Cesaro})) are
naturally connected to the Tomita-Takesaki modular operator (see subsection
\ref{examples}.3). In fact, this observation is common for p.d.f. with the
property ({\ref{kms}}). For those, we are going to produce a complete
classification in a subsequent paper.

\paragraph{\theparagraph.3.
The basic definition and the conjugate classes. } Let $\mathcal{H}$ be a
Hilbert space, $\mathcal{B}\left(\mathcal{H}\right)$ an algebra of all
bounded operators in $\mathcal{H}$, and $\mathcal{I}_{\mathcal{H}}$ the
identity operator in $\mathcal{H}$. We denote by
$\mathcal{U}\left(\mathcal{H}\right)$ the unitary subgroup in
$\mathcal{B}\left(\mathcal{H}\right)$. By a unitary representation of the
topological group $G$ we always mean a {\it continuous} homomorphism of $G$
into $\mathcal{U}\left(\mathcal{H}\right)$, where
$\mathcal{U}\left(\mathcal{H}\right)$ is equipped with the strong operator
topology.

\begin{Def}
A unitary representation $\pi:G\to\mathcal{U}\left(\mathcal{H}\right)$ of
$G$ is called a factor-representation if the $W^*-$algebra
$\pi(G)^{\prime\prime}$ generated by the operators $\pi(g)\;\left(g\in
G\right)$, is a factor.
\end{Def}
\begin{Def}
A unitary representation $\pi$ is called a factor-representation of finite
type if $\pi(G)^{\prime\prime}$ is a factor of type $II_1$.
\end{Def}

Let $\mathcal{M}$ be a factor of type $II_1$ and $\mathcal{M}$ a subalgebra
of $\mathcal{B}(\mathcal{H})$. If $\pi(G)\subset\mathcal{U}(\mathcal{M})=
\mathcal{M}\bigcap\mathcal{U}(\mathcal{H})$ and $tr_\mathcal{M}$ is the
unique normal, normalized $\left( tr(I)=1\right)$ trace on $\mathcal{M}$,
then it determines a {\it character} $\phi^{\mathcal{M}}_{\pi}$ of $G$ by
$\phi^{\mathcal{M}}_{\pi}(g) =tr_{\mathcal{M}}\left(\pi(g)\right)$.

\begin{Def}\label{trace}
A continuous function $\phi$ on $G$ is called a character if it satisfies
the following properties:
\begin{itemize}
\item [\it (a)] $\phi$ is central, that is, $\phi\left(g_1
g_2\right)=\phi\left(g_2 g_1\right)$ $\forall$ $g_1,g_2\in G$;

\item [\it (b)] $\phi$ is positive definite, that is, for all
$g_1,g_2,\ldots,g_n$ the matrix
$\left[\phi\left(g_jg_k^{-1}\right)\right]_{j,k=1}^n$ is non-negatively
definite;

\item [\it (c)] $\phi$ is normalized, that is, $\phi\left(e_G\right)=1$,
where $e_G$ is the unit of $G$.
\end{itemize}
\end{Def}

\begin{Def}\label{def4}
A character $\phi$ is called {\it indecomposable} if the group
representation corresponding to $\phi$ (according to the GNS construction)
is a factor-representation.
\end{Def}

In this paper we obtain a complete description of indecomposable characters
on $\Gamma\wr\mathfrak{S}_\infty$ in the case when $\Gamma$ is a separable
topological group.

First, let us describe the conjugacy classes in
$\Gamma\wr\mathfrak{S}_\infty$. Recall that the conjugacy classes in
$\mathfrak{S}_\infty$ are parameterized by {\it partitions} $\lambda$, that
is by unordered $\infty-$tuples $\lambda_1,\lambda_2,\ldots$ of natural
numbers such that there are finitely many elements $\lambda_k$ not equal to
$1$. Namely, $\lambda_1,\lambda_2,\ldots$ are the orders of cycles of a
permutation $s\in\mathfrak{S}_\infty$. Furthermore, an element
$\Gamma\wr\mathfrak{S}_\infty$ can be written as a product of an element of
$\mathfrak{S}_\infty$ and an element of $\Gamma^\infty_e$, and the
commutation rule between these two kinds of elements is as follows:
\begin{equation}\label{product}
s\cdot\gamma=s\cdot\left(\gamma_1,\gamma_2,\ldots
\right)=\left(\gamma_{s^{-1}(1)},\gamma_{s^{-1}(2)},\ldots\right)\cdot s,
\end{equation}
where $s\in\mathfrak{S}_\infty$,
$\gamma=\left(\gamma_1,\gamma_2,\ldots\right)\in\Gamma^\infty_e$. Let
$\mathbb{N}\diagup s$ be the set of orbits of $s$ on $\mathbb{N}$. Note
that for $p\in \mathbb{N}\diagup s$ the permutation $s_p$ given by
\begin{equation*}
s_p(k)=\left\{
\begin{array}{rl}
s(k)&\text{ if }k\in p
\\ k&\textit{otherwise}
\end{array}\right.,
\end{equation*}
is a cycle of order $|p|$, where $|p|$ stand for the cardinality of $p$.
For $\gamma=\left(\gamma_1,\gamma_2,\ldots\right)\in\Gamma^\infty_e$ define
the element $\gamma(p)=\left(\gamma_1(p),\gamma_2(p),\ldots \right)
\in\Gamma^\infty_e $ as follows:
\begin{equation}\label{color1}
\gamma_k(p)=\left\{
\begin{array}{rl}
\gamma_k&\text{ if }k\in p
\\ e&\textit{otherwise}.
\end{array}\right.
\end{equation}
Thus, using ({\ref{product}}), we have
\begin{equation}\label{decompositiontocycles}
s\cdot\gamma=\prod\limits_{p\in\mathbb{N}\diagup s}s_p\cdot\gamma(p).
\end{equation}
Denote by $\mathfrak{c}_{_G}(g)$ the conjugacy class of $g\in G$. Finally,
for $g=s\cdot\gamma\in\Gamma\wr\mathfrak{S}_\infty$ define the invariant
$\mathfrak{i}(g)$ given by unordered $\infty-$tuples of pairs
$\left\{\left(|p|,\mathfrak{c}_{_\Gamma}\left(\gamma_k\cdot\gamma_{s(k)}
\cdots\gamma_{s^{(l)}(k)}\cdots\gamma_{s^{(|p|-1)}(k)}\right)\right)
\right\}_{k\in p,\;p\in\mathbb{N}\diagup s}$, where $s^{(l)}$ is $l$-th
iteration of $s$. The following statement can be easily proved.

\begin{Prop}\label{conjugateclasses}
Let $g_1$ and $g_2$ be elements of $\Gamma\wr\mathfrak{S}_\infty$. Then
$\mathfrak{c}\left(g_1\right)=\mathfrak{c}\left(g_2\right)$ if and only if
$\mathfrak{i}\left(g_1\right)=\mathfrak{i}\left(g_2\right)$.
\end{Prop}

\paragraph{\theparagraph.4. The multiplicativity.}
The following claim gives a useful characterization of the class of
indecomposable characters:

\begin{Prop}\label{multiplicativity}
The following assumptions on a character $\phi$ of
$\Gamma\wr\mathfrak{S}_\infty$ are equivalent:
\begin{itemize}
\item [\it (a)] $\phi$ is indecomposable;

\item [\it (b)] $\phi(g)=\prod\limits_{p\in\mathbb{N}\diagup
s}\phi\left(s_p\cdot\gamma(p)\right)$ for any
$g=s\cdot\gamma=\prod\limits_{p\in\mathbb{N}\diagup s}s_p\cdot\gamma(p)$
(see {\ref{decompositiontocycles}}).
\end{itemize}
\end{Prop}

\begin{proof}
To prove the proposition, we consider the elements $g=s\cdot\gamma$ and
$g^{\prime}=s^{\prime}\cdot\gamma^{\prime}$ of
$\Gamma\wr\mathfrak{S}_\infty$ satisfying the following condition:
\begin{eqnarray*}
\left(\bigcup\limits_{p\in\left\{q\in\mathbb{N}\diagup
s\big|\,|q|>1\right\}}p\;\bigcup\left\{k|\gamma_k\neq
e\right\}\right)\subset\left(\bigcup\limits_{p\in\left\{q\in\mathbb{N}
\diagup s^\prime\big|\,|q|=1\right\}}p\;\bigcap\left\{k|\gamma_k^{\prime}=e
\right\}\right).
\end{eqnarray*}
Then there exists a sequence
$\left\{s_n\right\}_{n\in\mathbb{N}}\subset\mathfrak{S}_\infty$ such that
\begin{eqnarray}\label{assympt}
s_n\cdot g=g\cdot s_n\text{ and }s_n\cdot g^{\prime}s_n^{-1}\cdot h=h\cdot
s_n\cdot g^{\prime}s_n^{-1}\text{ for all }h\in \Gamma\wr\mathfrak{S}_n.
\end{eqnarray}
Suppose now that ({\it a}) holds. Using the GNS-construction, we produce
the representation $\pi_{\phi}$ of $\Gamma\wr\mathfrak{S}_\infty$ which
acts in a Hilbert space $\mathcal{H}_\phi$ with a cyclic vector $\xi_\phi$
such that
\begin{eqnarray*}
\phi(g)=\left(\pi_\phi\left(g\right)\xi_\phi,\xi_\phi\right).
\end{eqnarray*}
Let $A=w-\lim\limits_{n\to\infty}\pi_\phi\left(s_n\cdot
g^{\prime}s_n^{-1}\right)$ be a limit of the sequence
$\pi_\phi\left(s_n\cdot g^{\prime}s_n^{-1}\right)$ in the {\it weak
operator} topology. Using ({\ref{assympt}}), we deduce by Definition
{\ref{def4}} that $A=a\mathcal{I}$, where $\mathcal{I}$ is the identity
operator in $\mathcal{H}_\phi$ and $a$ a complex number. Therefore,
\begin{eqnarray*}
\phi\left(g\cdot g^\prime\right)=\lim\limits_{n\to\infty}\phi\left(g\cdot
s_n\cdot g^\prime\cdot s_n^{-1}\right)=\phi(g)\cdot\lim\limits_{n\to\infty}
\phi\left(s_n\cdot g^\prime\cdot
s_n^{-1}\right)=\phi(g)\cdot\phi\left(g^\prime\right).
\end{eqnarray*}
Thus {\it (b)} follows from {\it(a)}.

Conversely, suppose that {\it(b)} holds. For any subset $\mathcal{S}$ of
$\mathcal{B}\left(\mathcal{H}\right)$, define its commutant as follows:
\begin{eqnarray*}
\mathcal{S}^\prime=\left\{T\in\mathcal{B}\left(\mathcal{H}\right)\big|ST=TS
\text{ for all }S\in\mathcal{S}\right\}.
\end{eqnarray*}
If
$\pi_{\phi}\left(\Gamma\wr\mathfrak{S}_\infty\right)^\prime\bigcap\pi_{\phi}
\left(\Gamma\wr\mathfrak{S}_\infty\right)^{\prime\prime}=\mathcal{Z}$ is
larger than the scalars, then it contains a pair of orthogonal projections
$E$ and $F$ with the properties:
\begin{eqnarray}\label{0inequalities}
\phi\left(E\right)\ne 0,\phi\left(F\right)\ne 0\text{ and }E\cdot F=0.
\end{eqnarray}
By the von Neumann Double Commutant Theorem, for any $\varepsilon>0$ there
exist
$g_k^E,g_k^F\in\Gamma\wr\mathfrak{S}_n\subset\Gamma\wr\mathfrak{S}_\infty$
$\left(n<\infty\right)$ and complex numbers $c_k^E,c_k^F$
$\left(k=1,2,\ldots,N<\infty\right)$ such that
\begin{eqnarray}\label{inequalities}
\begin{split}
\left|\left|\sum\limits_{k=1}^Nc_k^E\pi_\phi\left(g_k^E\right)\xi_\phi
-E\xi_\phi\right|\right|<\varepsilon\phi(E),
\\ \left|\left|\sum\limits_{k=1}^N c_k^F\pi_\phi\left(g_k^F\right)
\xi_\phi-F\xi_\phi\right|\right|<\varepsilon\phi(F).
\end{split}
\end{eqnarray}
Consider the bijection
\begin{eqnarray*}
\tau(j)=\left\{
\begin{array}{rl}
j+n&\text{ if }j\le n,
\\ j-n&\text{ if }n<j\le 2n,
\\ j&\textit{otherwise}.
\end{array}\right.
\end{eqnarray*}
By Definition ({\ref{trace}}), use ({\ref{inequalities}}) to obtain
\begin{eqnarray}\label{1inequalities}
\begin{split}
\left|\left|\sum\limits_{k=1}^N c_k^E\pi_\phi\left(\tau
g_k^E\tau\right)\xi_\phi-E\xi_\phi\right|\right|<\varepsilon\phi(E),
\\ \left|\left|\sum\limits_{k=1}^Nc_k^F\pi_\phi\left(\tau g_k^F\tau\right)
\xi_\phi-F\xi_\phi\right|\right|<\varepsilon\phi(F).
\end{split}
\end{eqnarray}
Now, using {\it (b)}, ({\ref{0inequalities}}), ({\ref{inequalities}}) and
({\ref{1inequalities}}), we have
\begin{eqnarray*}
&&\varepsilon\sqrt{\phi(E)\phi(F)}\left(\varepsilon\sqrt{\phi(E)\phi(F)}+
\sqrt{\phi(E)}+\sqrt{\phi(F)}\right)
\\ &>&\left|\left(\sum\limits_{k=1}^Nc_k^E\pi_\phi\left(\tau g_k^E\tau
\right)\cdot \sum\limits_{k=1}^N c_k^F\pi_\phi \left( g_k^F
\right)\xi_\phi,\xi_\phi\right)\right|\\&=&\left|\left(\sum\limits_{k=1}^N
c_k^E\pi_\phi\left(\tau g_k^E\tau\right)\xi_\phi,\xi_\phi\right)\cdot
\left(\sum\limits_{k=1}^Nc_k^F\pi_\phi\left(\tau g_k^F\tau\right)
\xi_\phi,\xi_\phi \right)\right|
\\ &>&\phi(E)\phi(F)(\varepsilon+1)^2.
\end{eqnarray*}
Hence
\begin{eqnarray*}
\varepsilon>\left[\frac{1-\sqrt{\phi(F)}}{\sqrt{\phi(F)}}+
\frac{1-\sqrt{\phi(E)}}{\sqrt{\phi(E)}}\right]^{-1}.
\end{eqnarray*}
Then, comparing this to ({\ref{0inequalities}}), we get a contradiction.
\end{proof}

\paragraph{\theparagraph.5. The main result.}
In {\cite{Thoma}}, E. Thoma obtained the following remarkable description
of all {\it indecomposable} characters of $\mathfrak{S}_\infty$. The
characters of $\mathfrak{S}_\infty$ are labeled by pairs of non-increasing
positive sequences of numbers $\left\{\alpha_k\right\}$,
$\left\{\beta_k\right\}$ $\left(k\in\mathbb{N}\right)$ such that
\begin{eqnarray}\label{cond}
\sum\limits_{k=1}^{\infty}\alpha_k+\sum\limits_{k=1}^{\infty}\beta_k\le 1.
\end{eqnarray}
The value of the corresponding character on a permutation with a single
cycle of length $l$ is
\begin{eqnarray*}\label{color2}
\sum\limits_{k=1}^{\infty}\alpha_k^l
+(-1)^{l-1}\sum\limits_{k=1}^{\infty}\beta_k^l.
\end{eqnarray*}
Its value on a permutation with several disjoint cycles equals to the
product of values on each cycle.

Here is our main result.
\begin{Th}\label{mainth}
Let $\phi$ be an indecomposable character of
$\Gamma\wr\mathfrak{S}_\infty$. Then there exist a representation
$\varrho^0$ of $\Gamma$ of {\it finite} type, two non-increasing sequences
of positive numbers $\left\{\alpha_k\right\}$, $\left\{\beta_k\right\}$
$\left(k\in\mathbb{N}\right)$, and two sequences
$\left\{\varrho^{\alpha_k}\right\}$ and $\left\{\varrho^{\beta_k}\right\}$
of finite-dimensional irreducible representations of $\Gamma$ such that for
$g=s\cdot\gamma\in \Gamma\wr\mathfrak{S}_\infty$ (see ({\ref{product}}) --
({\ref{color1}})) one has
\begin{eqnarray}
\begin{split}
&\phi(g)=\prod\limits_{p\in\mathbb{N}\diagup
s}\Bigg\{\delta_p\cdot\left[1-\sum\limits_{k\in\mathbb{N}}
\left(\alpha_k\cdot\dim\varrho^{\alpha_k}+
\beta_k\cdot\dim\varrho^{\beta_k}\right)\right]\cdot
\prod\limits_{j\in p}\;tr_0\left( \gamma_j \right)
\\ &+\sum\limits_{k=1}^{\infty}\Bigg[\alpha_k^{|p|}\cdot
Tr_{\alpha_k}\left(\tilde{\gamma}(p)\right)+(-1)^{|p|-1}\beta_k^{|p|}\cdot
Tr_{\beta_k}\left(\tilde{\gamma}(p)\right)\Bigg]\Bigg\},
\end{split}
\end{eqnarray}
where
$\tilde{\gamma}(p)=\gamma_k\cdot\gamma_{s(k)}\cdots\gamma_{s^{(l)}(k)}\cdots
\gamma_{s^{(|p|-1)}(k)}$ $(k\in p)$; $Tr_{\alpha_k}$, $Tr_{\beta_k}$ are
characters corresponding to the representations $\varrho^{\alpha_k}$,
$\varrho^{\beta_k}$, $tr_0$ is the normalized character of the
representation $\varrho^0$;
\begin{eqnarray*}
\delta_{p}=\left\{
\begin{array}{rl}
1&\text{ if }|p|=1,
\\ 0&\text{ if }|p|>1.
\end{array}\right.
\end{eqnarray*}
and
$\sum\limits_k\left(\alpha_k\cdot\dim\varrho^{\alpha_k}+\beta_k\cdot\dim
\varrho^{\beta_k}\right)\le 1$.
\end{Th}
Now we formulate the main result in the case when $\Gamma$ is a {\it
locally compact} abelian group.

Let $\overset{\land}{\Gamma}$ stand for the {\it dual} group of $\Gamma$.

\begin{Th}
Let $\phi$ be an indecomposable character of
$\Gamma\wr\mathfrak{S}_\infty$. There exist a probability measure $\mu$ on
$\overset{\land}{\Gamma}$, two non-increasing sequences of positive numbers
$\left\{\alpha_k\right\}$, $\left\{\beta_k\right\}$
$\left(k\in\mathbb{N}\right)$, two sequences
$\left\{\overset{\land}{\alpha}_k\right\}$ and
$\left\{\overset{\land}{\beta}_k\right\}$ of elements of
$\overset{\land}{\Gamma}$, such that for
$g=s\cdot\gamma\in\Gamma\wr\mathfrak{S}_\infty$ (see ({\ref{product}}) --
({\ref{color1}}))
\begin{eqnarray}
\begin{split}
&\phi(g)=\prod\limits_{p\in\mathbb{N}\diagup
s}\Bigg\{\delta_p\;\left[1-\sum\limits_{k\in\mathbb{N}}\left(\alpha_k+
\beta_k\right)\right]\cdot\int\limits_{\overset{\land}{\gamma}\in
\overset{\land}{\Gamma}}\overset{\land}{\gamma}
\left(\prod\limits_{j\in\mathbb{N}}\gamma_j(p)\right)\,d\,\mu
\\ &+\sum\limits_{k=1}^{\infty}\Bigg[\alpha_k^{|p|}\cdot
\overset{\land}{\alpha}_k\left(\prod\limits_{j\in\mathbb{N}}\gamma_j(p)
\right)+(-1)^{|p|-1}\beta_k^{|p|}\cdot\overset{\land}{\beta}_k
\left(\prod\limits_{j\in\mathbb{N}}\gamma_j(p)\right)\Bigg]\Bigg\},
\end{split}
\end{eqnarray}
where $\left\{\alpha_k\right\}$ and $\left\{\beta_k\right\}$ satisfy
({\ref{cond}}).
\end{Th}

\section{Realizations of $II_1-$factor-representations.}\label{two} A
complete family of $II_1-$factor-representations of $G\wr\Gamma$ can be
constructed using the Vershik-Kerov {\cite{VK1}} or Olshanski {\cite{O2}}
realizations, found for the $ II_1-$factor-representations of the infinite
symmetric group $\mathfrak{S}_\infty$. We follow the approach developed by
Olshanski as it leads to less spadework.

\paragraph{\theparagraph.1. A construction of representations.}

Let $\left\{\alpha_k\right\}_{k\in\mathbb{N}}$,
$\left\{\beta_k\right\}_{k\in\mathbb{N}}$ be two finite or infinite sets of
numbers from $(0,1)$ and suppose that $\varrho^{\alpha_k}$ and
$\varrho^{\beta_k}$ are unitary irreducible finite dimensional
representations of $\Gamma$ that act in the Hilbert spaces
$\mathcal{H}^{\alpha_k}$ and $\mathcal{H}^{\beta_k}$ respectively. Assume
that
\begin{eqnarray*}
\sum\limits_k\alpha_k\cdot\dim
\varrho^{\alpha_k}+\sum\limits_k\beta_k\cdot\varrho^{\beta_k}\leq 1.
\end{eqnarray*}
Set
\begin{eqnarray*}
\delta=1-\sum\limits_k\alpha_k\cdot\dim\varrho^{\alpha_k}-
\sum\limits_k\beta_k\cdot\varrho^{\beta_k}.
\end{eqnarray*}
Let $\mathcal{H}^0$ stand for the (Hilbert) space of a unitary
representation $\varrho^0$ of $\Gamma$ of finite type. We may assume
without loss of generality that there exists a {\it cyclic} and {\it
separating} unit vector $\xi^{(0)}$ for the pair
$\left(\varrho^0\left(\Gamma\right),\mathcal{H}^0\right)$. To rephrase,
$\left[\varrho^0(\Gamma)\xi^{(0)}\right]=
\left[\varrho^0(\Gamma)^\prime\xi^{(0)}\right]=\mathcal{H}^0$, where
$\left[\varrho^0(\Gamma)\xi^{(0)}\right]$ is the subspace generated by
$\varrho^0(\Gamma)\xi^{(0)}$. Furthermore, the formula
$tr^0(\gamma)=\left(\varrho^0(\gamma)\xi^{(0)},\xi^{(0)}
\right)_{\mathcal{H}^0}$ determines a character on $\Gamma$. Thus, we
associate a unitary representation $\left(\varrho^0\right)^{(2)}$ of
$\Gamma\times\Gamma$ to the representation $\varrho^0$. Namely,
$\left(\varrho^0\right)^{(2)}$ is defined as follows:
\begin{eqnarray*}
\left(\varrho^0\right)^{(2)}\left((\gamma_1,\gamma_2)\right)
\left(\varrho^0(\gamma)\xi^{(0)}\right)=\varrho^0\left(\gamma_1\right)
\varrho^0(\gamma)\varrho^0\left(\gamma_2^{-1}\right)\xi^{(0)}.
\end{eqnarray*}
Denote by $\left(\varrho^{0k},\mathcal{H}^{0k},\xi^{(0k)}\right)$ the k-th
copy of the triplet $\left(\varrho^0,\mathcal{H}^0,\xi^{(0)}\right)$.

Let $\left\{\mathrm{e}_j^{(\alpha_k)}\right\}_{1\le
j\le\dim\mathcal{H}^{\alpha_k}}$ be an orthonormal basis in
$\mathcal{H}^{\alpha_k}$. Define the matrix elements of
$\overline{\varrho}^{\alpha_k}$ as follows:
\begin{eqnarray*}
\overline{\varrho}^{\alpha_k}_{j\,l}(\gamma)=
\overline{\left(\varrho^{\alpha_k}(\gamma)\mathrm{e}_j^{(\alpha_k)},\;
\mathrm{e}_l^{(\alpha_k)}\right)},
\end{eqnarray*}
where bar denotes the complex conjugation.

Let
\begin{eqnarray*}
\mathbf{H}=
\left(\underset{k}{\oplus}\mathcal{H}^{\alpha_k}\underset{k}{\oplus}
\mathcal{H}^{\beta_k}\right)\otimes\left(\underset{k}{\oplus}
\mathcal{H}^{\alpha_k}\underset{k}{\oplus}\mathcal{H}^{\beta_k}\right)
\underset{k}{\oplus}\mathcal{H}^{0k}
\end{eqnarray*}
and
\begin{eqnarray*}
\eta^{(m)}=\sum\limits_k\alpha_k\left(\sum\limits_{j}
\mathrm{e}_j^{(\alpha_k)}\otimes\mathrm{e}_j^{(\alpha_k)}\right)+
\sum\limits_k\beta_k\left(\sum\limits_{j}\mathrm{e}_j^{(\beta_k)}\otimes
\mathrm{e}_j^{(\beta_k)}\right)+\sqrt{\delta}\xi^{(0m)}.
\end{eqnarray*}
Define the unitary representations $\varrho$ and $\bar{\varrho}\,$ of
$\Gamma$ in $\mathbf{H}$ as follows
\begin{eqnarray}
\begin{split}
&\varrho=\left(\underset{k}{\oplus}\varrho^{\alpha_k}\underset{k}{\oplus}
\varrho^{\beta_k}\right)\otimes\left(\underset{k}{\oplus}I
\underset{k}{\oplus}I\right)\underset{k}{\oplus}\varrho^{0k}
\\ &\bar{\varrho}=\left(\underset{k}{\oplus}I\underset{k}{\oplus}I\right)
\otimes\left(\underset{k}{\oplus}\bar{\varrho}^{\alpha_k}
\underset{k}{\oplus}\bar{\varrho}^{\beta_k}\right)
\underset{k}{\oplus}\bar{\varrho}^{0k},\;\;\textit{ where }
\;\;\bar{\varrho}^{0k}(\gamma)=\left(\varrho^0\right)^{(2)}
\left(\left(e_\Gamma,\gamma\right)\right).
\end{split}
\end{eqnarray}

We identify $\mathcal{H}^{\alpha_k}\otimes \mathcal{H}^{\alpha_k}$,
$\mathcal{H}^{\beta_k}\otimes\mathcal{H}^{\beta_k}$, and $\mathcal{H}^{0k}$
to their images with respect to their natural embeddings to $\mathbf{H}$.
Denote by $\mathbf{H}^{{m}}$ the $m-$th copy of the Hilbert space
$\mathbf{H}$ and consider the infinite tensor product
\begin{eqnarray*}
\overset{\smallsmile}{\mathbf{H}}=\bigotimes\limits_m
\left(\mathbf{H}^{m},\eta^{(m)}\right).
\end{eqnarray*}
It is convenient to represent $\overset{\smallsmile}{\mathbf{H}}$ as the
closure of linear span of vectors of the form
\begin{eqnarray*}
\zeta_1\otimes\zeta_2\otimes\cdots\otimes\zeta_{m-1}\otimes\eta^{(m)}
\otimes\eta^{(m+1)}\otimes\cdots,\text{ where }\zeta_j\text{ is any vector
from }\mathbf{H}^{j}.
\end{eqnarray*}
Now fix the orthonormal basis
\begin{eqnarray*}
\mathfrak{B}=\left\{\mathrm{e}_j^{(r)}\otimes\mathrm{e}_l^{(s)}\in
\left(\underset{k}{\oplus}\mathcal{H}^{\alpha_k}\underset{k}{\oplus}
\mathcal{H}^{\beta_k}\right)\otimes\left(\underset{k}{\oplus}
\mathcal{H}^{\alpha_k}\underset{k}{\oplus}\mathcal{H}^{\beta_k}\right),
\mathrm{e}_j\in\underset{k}{\oplus}\mathcal{H}^{0k}\right\}
\end{eqnarray*}
in $\mathbf{H}$ and assume below $\zeta_j\in\mathfrak{B}$. By the vector
$\overset{\smallsmile}{\zeta}=\zeta_1\otimes\zeta_2\otimes\cdots\otimes
\zeta_{m-1}\otimes\cdots$ we build a sequence
$\mathfrak{j}\left(\overset{\smallsmile}{\zeta}\right)=\left\{j_1<j_2<\cdots
\right\}$ such, that
\begin{eqnarray*}
\zeta_{j_{\,l}}=\mathrm{e}_m^{(\beta_k)}\otimes\mathrm{f}\textit{ for some
}\beta_k\textit{ and }m.
\end{eqnarray*}
Define for $\mathfrak{s}\in\mathfrak{S}_\infty$ a vector
$\mathfrak{s}\left({\overset{\smallsmile}{\zeta}}\right)=
\vartheta_1\otimes\vartheta_2\otimes\cdots\otimes\vartheta_{m-1}\otimes
\cdots$ as follows:
\begin{eqnarray*}
&\vartheta_k=\left\{
\begin{array}{rl}
\mathrm{e}_{i_k}^{(r)}\otimes\mathrm{f},&\textit{ if }\zeta_{s^{-1}(k)}=
\mathrm{e}_{i_k}^{(r)}\otimes\mathrm{e}\;\textit{ and }\;
\zeta_k=\mathrm{e}_j^{(s)}\otimes\mathrm{f}
\\ \mathrm{e}_{i_k}^{(r)}\otimes\mathrm{e},&\textit{ if }\zeta_{s^{-1}(k)}=
\mathrm{e}_{i_k}^{(r)}\otimes\mathrm{e}\;\textit{ and }\;
\zeta_k\in\underset{l}{\oplus}\mathcal{H}^{0l}
\\ \zeta_{s^{-1}(k)},&\textit{ if }\zeta_{s^{-1}(k)}\in\underset{l}{\oplus}
\mathcal{H}^{0l}.
\end{array}\right.
\end{eqnarray*}
For any $j_l\in\mathfrak{j}\left(\overset{\smallsmile}{\zeta}\right)$ such
that $\zeta_{j_{\,l}}=
\mathrm{e}_{m\left(j_l\right)}^{(\beta_k)}\otimes\mathrm{f}$ there exists
$j_l^\mathfrak{s}\in\mathfrak{j}\left(\mathfrak{s}
\left(\overset{\smallsmile}{\zeta}\right)\right)$ with the property
\begin{eqnarray*}
\vartheta_{j_l^\mathfrak{s}}=\mathrm{e}_{m\left(
j_l\right)}^{(\beta_k)}\otimes\mathrm{g}
\end{eqnarray*}
Let $\mathfrak{t}$ be a permutation of the set $\left\{
j_1^\mathfrak{s},j_2^\mathfrak{s}, \ldots \right\}$ for which $\mathfrak{t}
\left(j_1^\mathfrak{s}\right)< \mathfrak{t}
\left(j_2^\mathfrak{s}\right)<\ldots $. Finally, set
$\psi\left(\mathfrak{s},
\overset{\smallsmile}{\zeta}\right)=sgn(\mathfrak{t})$. The corresponding
representation $\pi$ of $\mathfrak{S}_\infty\wr\Gamma$ can be realized in
$\overset{\smallsmile}{\mathbf{H}}$ as follows:
\begin{eqnarray}
\begin{split}
&\pi(\gamma)\left(\zeta_1\otimes\zeta_2\otimes\cdots\otimes\zeta_{m-1}
\otimes\eta^{(m)}\otimes\cdots\right)
\\ =&\varrho\left(\gamma_1\right)\zeta_1\otimes\varrho
\left(\gamma_2\right)\zeta_2\otimes\cdots\otimes\varrho\left(\gamma_{m-1}
\right)\zeta_{m-1}\otimes\varrho\left(\gamma_m\right)\eta^{(m)}\otimes\cdots
\\ &\textit{ and for }\mathfrak{s}\in\mathfrak{S}_\infty\;\;\;
\pi(\mathfrak{s})\left(\zeta_1\otimes\zeta_2\otimes\cdots\otimes
\zeta_{m-1}\otimes\cdots\right)=
\psi\left(\mathfrak{s},\overset{\smallsmile}{\zeta}\right)\mathfrak{s}
\left({\overset{\smallsmile}{\zeta}}\right).
\end{split}
\end{eqnarray}

\paragraph{\theparagraph.2. The character's formula.}
Set $\overset{\smallsmile}{\eta}=\bigotimes\limits_m\eta^{(m)}$. Assume
that $\mathfrak{s}$ is the cycle $(1\to 2\to 3\to\cdots\to k-1\to k)$. Let
$\gamma=\left(\gamma_1,\gamma_2,\ldots,\gamma_k,e_\Gamma,e_\Gamma,\ldots
\right)$. Routine calculations provide that
\begin{eqnarray}\label{realchar}
\left(\pi\left(\mathfrak{s}\gamma\right)
\overset{\smallsmile}{\eta},\overset{\smallsmile}{\eta}\right)=
\sum\limits_j{\alpha_j^k}Tr\left(\varrho^{\alpha_j}
\left(\gamma_1\gamma_2\cdots\gamma_k\right)\right)+\sum\limits_jTr
\left(\varrho^{\beta_j}\left(\gamma_1\gamma_2\cdots \gamma_k\right)\right),
\end{eqnarray}
where $Tr\left(\varrho^r(\gamma)\right)=\sum\limits_{j=1}^{\dim\,\varrho^r}
\varrho_{j\,j}^r(\gamma)$ and $k>1$.

It is obvious, that
\begin{eqnarray*}
\left(\pi\left(\gamma\right)\overset{\smallsmile}{\eta},
\overset{\smallsmile}{\eta}\right)=
\prod\limits_{j=1}^k\left(\sum\limits_i{\alpha_i}Tr\left(\varrho^{\alpha_i}
\left(\gamma_j\right)\right)+\sum\limits_i
{\beta_i}Tr\left(\varrho^{\beta_i}\left(\gamma_j\right)\right)+
\left(\varrho^0\left(\gamma_j\right)\xi^{(0)},\xi^{(0)}\right)\right).
\end{eqnarray*}
Since $tr^0$ is a character on $\Gamma$, one can use ({\ref{realchar}}) and
the multiplicativity property (see Proposition {\ref{multiplicativity}}) to
obtain the following
\begin{Co}
Let $\chi(g)=\left(\pi\left(g\right)\overset{\smallsmile}{\eta},
\overset{\smallsmile}{\eta}\right)$. Then $\chi$ is an indecomposable
character on $\mathfrak{S}_\infty\wr\Gamma$.
\end{Co}

\section{Another examples.}\label{examples} In this section we construct
examples of {\it infinite} type representations of
$\mathfrak{S}_\infty\wr\mathbb{Z}_2$. The corresponding positive definite
functions are not characters. On the other hand they satisfy the following
condition:
\begin{eqnarray*}
\varphi(sg)=\varphi(gs)\textit{ for all }g\in
G=\mathfrak{S}_\infty\wr\Gamma\textit{ and }s\in\mathfrak{S}_\infty.
\end{eqnarray*}
In the generic case the representation $\pi_\varphi$ built by
GNS-construction from $\varphi$ is of type $I\!I\!I$. Furthermore, the
state $\varphi$ on the $W^*-$algebra
$\pi_\varphi\left(G\right)^{\prime\prime}$ is exact. These properties allow
one to construct the Tomita-Takesaki modular operator $\Delta_\varphi$.
Surprisingly, $\Delta_\varphi$ is naturally related to the Okounkov
operator $\mathcal{O}_k$ (see (\ref{Cesaro})), which is an important in the
representation theory of symmetrical group (see \cite{Ok1}, \cite{Ok2}).

\paragraph{\theparagraph.1. A construction.} Let
$X_i=\mathbb{Z}_2\times\mathbb{Z}_2=\{0,1\}\times\{0,1\}$. Define a
probability measure $\nu_i$ on $X_i$ by $\nu_i((k,l))=p_{kl}$. Let
$\left(X,\mu\right)=\prod\limits_i\left(X_i,\nu_i\right)$ and
$x=\left(x_i\right)\in X$, where $x_i=\left(x_i^{(0)},x_i^{(1)}\right)\in
X_i$, $x_i^{(k)}\in\{0,1\}$. Define an action $\mathfrak{a}$ of
$g=\left(s_0,s_1\right)\in\mathfrak{S}_\infty\times\mathfrak{S}_\infty$ on
$\left(X,\mu\right)$ as follows:
\begin{eqnarray*}
\left(\mathfrak{a}_g(x)\right)_i^{(k)}=x_{s_k(i)}^{(k)}\;\;\;(k=0,1).
\end{eqnarray*}
\label{pegrem}
\begin{Rem}\label{Rem1}
The measure $\mu$ is
$\mathfrak{G}_\infty\times\mathfrak{G}_\infty-$quasiinvariant if and only
if $p_{ij}\ne 0$ for all $i,j=0,1$.
\end{Rem}
We are about to construct a unitary representation $\pi_\mu$ of $G\times G$
in $L^2\left(X,\mu\right)$. With $\varsigma\in L^2\left(X,\mu\right)$ set
up
\begin{eqnarray}
\begin{split}
\left(\pi_\mu\left(\left(s_0,s_1\right)\right)\varsigma\right)(x)=
\left(\frac{d\,\mu\left(\mathfrak{a}_g(x)\right)}{d\,\mu\left(x\right)}
\right)^{\frac{1}{2}}\varsigma\left(\mathfrak{a}_g(x)\right),
\\ \left(\pi_\mu\left(\left(\gamma^{(0)},\gamma^{(1)}\right)\right)\varsigma
\right)(x)=(-1)^{\left(\sum\limits_{i,k}\gamma_i^{(k)}x_i^{(k)}\right)}
\varsigma (x),
\end{split}
\end{eqnarray}
where $\gamma^{(0)}=\left(\gamma_i^{(0)}\right)\in\mathbb{Z}_2^\infty$,
$\gamma^{(1)}=\left(\gamma_i^{(1)}\right)\in\mathbb{Z}_2^\infty$, and
$\left(\gamma^{(0)},\gamma^{(1)}\right)\in\mathbb{Z}_2^\infty\times
\mathbb{Z}_2^\infty$. Let $\pi_\mu^{(0)}(g)=\pi_\mu\left(\left(g,e_G
\right)\right)$ and
$\pi_\mu^{(1)}(g)=\pi_\mu\left(\left(e_G,g\right)\right)$.
\begin{Prop}
$\pi_\mu$ is irreducible. Hence, $\pi_\mu^{(0)}$ and $\pi_\mu^{(1)}$ are
factor-representations of $\mathfrak{S}_\infty\wr\Gamma$.
\end{Prop}
\begin{proof}
Obvious.
\end{proof}

\paragraph{\theparagraph.2. A cyclic separating vector.}
Let $\mathbb{I}$ be an element of $L^2\left( X,\mu \right)$ given by the
function identically equal to the unit.
\begin{Th} \label{cyclicth}
If $\det\left[p_{i\,j}\right]\ne 0$, then $\mathbb{I}$ is a cyclic
separating vector for $ \pi_\mu^{(0)}(G)^{\prime\prime}$ and
$\pi_\mu^{(1)}(G)^{\prime\prime}$. That is,
\begin{eqnarray*}
\left[\pi_\mu^{(0)}(G)^{\prime\prime}\mathbb{I}\right]=
\left[\pi_\mu^{(1)}(G)^{\prime\prime}\mathbb{I}\right]=
L^2\left(X,\mu\right).
\end{eqnarray*}
\end{Th}
\begin{proof}
Let $(k,l)$ be a transposition from $\mathfrak{S}_\infty$. First notice
that the operator
\begin{eqnarray*}
\mathcal{O}_k^{(j)}=\lim\limits_{n\to\infty}\frac{1}{n}\sum\limits_{l=1}^n
\pi_\mu^{(j)}\left((k,l)\right)\;\;(\textit{ see ({\ref{Cesaro}})})
\end{eqnarray*}
belongs to $\pi_\mu^{(j)}(G)^{\prime\prime}$ $(j=0,1)$. Since
\begin{eqnarray*}
\left(L^2\left(X,\mu\right),\mathbb{I}\right)=\bigotimes
\limits_{i=1}^\infty\left(L^2\left( X_i,\nu_i\right),\mathbb{I}\right)
\end{eqnarray*}
one can apply the law of large numbers to deduce that
\begin{eqnarray*}
\mathcal{O}_i^{(j)}=I\otimes I\otimes\ldots\otimes
\underset{i-\textit{th}}{\mathcal{O}_i^{(j,i)}}\otimes I\otimes\ldots.
\end{eqnarray*}
Furthermore, if $\chi_{kl}^{(i)}$ is the indicator of the point $(k,l)\in
X_i= \mathbb{Z}_2\times\mathbb{Z}_2$, the matrices of
$\mathcal{O}_i^{(0,i)}$ and $\mathcal{O}_i^{(1,i)}$ in the orthonormal
basis $\left\{ \mathrm{e}_{kl}^{(i)}=\frac{\chi_{kl}^{(i)}}{\sqrt{p_{kl}}}
\right\}_{k,l=0,1}$ are as follows:
\begin{eqnarray}\label{matrices}
\begin{split}
\mathcal{O}_i^{(0,i)}\leftrightarrow\left[
\begin{smallmatrix}
{p_{00}+p_{01}}&{0}&{\sqrt{p_{00}
p_{10}}+\sqrt{p_{01}p_{11}}}&{0}
\\ 0&p_{00}+p_{01}&0&\sqrt{p_{00}p_{10}}+\sqrt{p_{01}p_{11}}
\\ \sqrt{p_{00}p_{10}}+\sqrt{p_{01}p_{11}}&0&p_{10}+p_{11}&0
\\ 0&\sqrt{p_{00}p_{10}}+\sqrt{p_{01}p_{11}}&0&p_{10}+p_{11}
\end{smallmatrix}
\right],
\\ \mathcal{O}_i^{(1,i)}\leftrightarrow\left[
\begin{smallmatrix}
p_{00}+p_{10}&\sqrt{p_{00}p_{01}}+\sqrt{p_{10}p_{11}}&0&0
\\ \sqrt{p_{00}p_{01}}+\sqrt{p_{10}p_{11}}&p_{01}+p_{11}&0&0
\\ 0&0&p_{00}+p_{10}&\sqrt{p_{00}p_{01}}+\sqrt{p_{10}p_{11}}
\\ 0&0&\sqrt{p_{00}p_{01}}+\sqrt{p_{10}p_{11}}&p_{01}+p_{11}
\end{smallmatrix}
\right].
\end{split}
\end{eqnarray}
By the construction,
\begin{eqnarray*}
\pi_\mu^{(k)}\left(\gamma^{(k)}\right)=\bigotimes\limits_{i=1}^\infty
\pi_\mu^{(k,i)}\left(\gamma_i^{(k)}\right),
\end{eqnarray*}
where $ \pi_\mu^{(0,i)}\left(\gamma_i^{(0)}\right)$ and $\pi_\mu^{(1,i)}
\left(\gamma_i^{(1)}\right)$ are determined by the matrices
\begin{eqnarray}\label{matricesOk1}
\left[
\begin{matrix}
1&0&0&0
\\ 0&1&0&0
\\ 0&0&(-1)^{\gamma_i^{(0)}}&0
\\ 0&0&0&(-1)^{\gamma_i^{(0)}}
\end{matrix}
\right]\textit{ and }\left[
\begin{matrix}
1&0&0&0
\\ 0&(-1)^{\gamma_i^{(1)}}&0&0
\\ 0&0&1&0
\\ 0&0&0&(-1)^{\gamma_i^{(1)}}
\end{matrix}
\right].
\end{eqnarray}
Use the map
\begin{eqnarray}\label{isom}
\mathfrak{I}_i:\;\sum\limits_{m,n=0,1}a_{mn}\mathrm{e}_{mn}^{(i)}\to \left[
\begin{matrix}
a_{00}&a_{01}
\\ a_{10}&a_{11}
\end{matrix}\right],
\end{eqnarray}
to identify $L^2\left(X_i,\nu_i\right)$ to the full matrix algebra
$M_2\left(\mathbb{C}\right)$, so that
\begin{eqnarray*}
\mathfrak{I}_i(\mathbb{I})=\left[
\begin{matrix}
\sqrt{p_{00}}&\sqrt{p_{01}}
\\ \sqrt{p_{10}}&\sqrt{p_{11}}
\end{matrix}\right].
\end{eqnarray*}
Equip $M_2\left(\mathbb{C}\right)$ with the Hermitian form
\begin{eqnarray*}
\left<a,b\right>_i=Tr\left( b^*a\right),
\end{eqnarray*}
then $\mathfrak{I}_i$ is a unitary and
$\mathfrak{I}_iL^2\left(X_i,\nu_i\right)=M_2\left(\mathbb{C}\right)$. Now
as an elementary consequence of {(\ref{matrices})} and
{(\ref{matricesOk1})} one has:
\begin{eqnarray}\label{lr}
\begin{split}
&\mathfrak{I}_i\mathcal{O}_i^{(0,i)}\mathfrak{I}_i^{-1}a=\left[
\begin{matrix}
p_{00}+p_{01}&\sqrt{p_{00}p_{10}}+\sqrt{p_{01}p_{11}}
\\ \sqrt{p_{00}p_{10}}+\sqrt{p_{01}p_{11}}&p_{10}+p_{11}
\end{matrix}\right]a,
\\ & \mathfrak{I}_i\mathcal{O}_i^{(1,i)}\mathfrak{I}_i^{-1}a=a\left[
\begin{matrix}
p_{00}+p_{10}&\sqrt{p_{00}p_{01}}+\sqrt{p_{10}p_{11}}
\\ \sqrt{p_{00}p_{01}}+\sqrt{p_{10}p_{11}}&p_{01}+p_{11}
\end{matrix}\right],
\\ &\mathfrak{I}_i \pi_\mu^{(0,i)}\left(\gamma_i^{(0)}\right)
\mathfrak{I}_i^{-1}a=\left[
\begin{matrix}
1&0
\\ 0&(-1)^{\gamma_i^{(0)}}
\end{matrix}
\right]a,
\\ &\mathfrak{I}_i\pi_\mu^{(1,i)}\left(\gamma_i^{(1)}\right)
\mathfrak{I}_i^{-1}a=a\left[\begin{matrix}
1&0
\\ 0&(-1)^{\gamma_i^{(1)}}
\end{matrix}
\right],\textit{ where }a\in M_2(\mathbb{C}).
\end{split}
\end{eqnarray}
Thus, in view of Remark {\ref{Rem1}} (see p. \pageref{pegrem}), the algebra
$\mathfrak{M}_i^k$ generated by the operators
$\mathfrak{I}_i\mathcal{O}_i^{(k,i)}\mathfrak{I}_i^{-1}$ and
$\mathfrak{I}_i\pi_\mu^{(0,i)}\left(\gamma_i^{(k)}\right)
\mathfrak{I}_i^{-1}$ is just $M_2(\mathbb{C})$. Since
$\det\left(\mathfrak{I}_i\left(\mathbb{I}\right)\right)\ne 0$, one has
finally
$\mathfrak{M}_i^0\mathfrak{I}_i\left(\mathbb{I}\right)=\mathfrak{M}_i^1
\mathfrak{I}_i\left(\mathbb{I}\right)=M_2(\mathbb{C})$.
\end{proof}

\paragraph{\theparagraph.3. The modular operator.}
Consider the Hilbert space
$\mathfrak{H}=\bigotimes\limits_{i=1}^\infty\left(M_2(\mathbb{C}),\left<
\;\;\right>_i,\mathfrak{I}_i(\mathbb{I})\right)$. It is convenient to
represent $\mathfrak{H}$ as a closure of the linear span of the vectors
$a_1\otimes a_2\otimes\ldots\otimes
a_i\otimes\mathfrak{I}_{i+1}(\mathbb{I})\otimes
\mathfrak{I}_{i+2}(\mathbb{I})\ldots$, where $a_i\in M_2(\mathbb{C})$. If
$\mathfrak{I}=\bigotimes\limits_{i=1}^\infty\mathfrak{I}_i$, one has by
Theorem {\ref{cyclicth}}
\begin{eqnarray*}
\mathfrak{I}L^2\left(X,\mu\right)=\mathfrak{H}.
\end{eqnarray*}
Let $\mathcal{L}(\mathfrak{H})$ and $\mathcal{R}(\mathfrak{H})$ be the
$W^*-$algebras generated in $\mathfrak{H}$ by the operators of left and
right multiplication by elements of the form
\begin{eqnarray*}
a_1\otimes a_2\otimes\ldots\otimes a_i\otimes I_2\otimes
I_2\otimes\ldots,\textit{ where }a_i\in M_2({\mathbb{C}}),\;\;\;I_2=\left[
\begin{matrix}
1&0
\\ 0&1 \end{matrix}\right].
\end{eqnarray*}

\begin{Prop}
$\pi_\mu^{(0)}(G)^{\prime\prime}=\mathfrak{I}^{-1}
\mathcal{L}(\mathfrak{H})\mathfrak{I}$ and
$\pi_\mu^{(1)}(G)^{\prime\prime}=\mathfrak{I}^{-1}
\mathcal{R}(\mathfrak{H})\mathfrak{I}$.
\end{Prop}
\begin{proof}
Let $\mathfrak{A}_n^{(j)}$ stand for the $W^*-$algebra generated by the
operators $\left\{\mathcal{O}_i^{(j)}\right\}_{i=1}^n$ and
$\left\{\pi_\mu^{(j)}\left(\Gamma^n\right)\right\}$ $(j=0,1)$. In view of
{(\ref{lr})}, $\mathfrak{A}_n^{(j)}$ is isomorphic
$\bigotimes\limits_{i=1}^n M_2(\mathbb{C})$. Therefore,
$\pi_\mu^{(j)}\left(\mathfrak{S}_n\right)\subset \mathfrak{A}_n^{(j)}$.
Finally, use {(\ref{lr})} deduce
$\mathfrak{A}_n^{(0)}\subset\mathcal{L}(\mathfrak{H})$ and
$\mathfrak{A}_n^{(1)}\subset\mathcal{R}(\mathfrak{H})$.
\end{proof}

\medskip

Let $\xi=\mathfrak{I}_{1}(\mathbb{I})\otimes
\mathfrak{I}_{2}(\mathbb{I})\otimes\ldots\otimes
\mathfrak{I}_{i+2}(\mathbb{I})\otimes\ldots$. Since the vector $\xi$ is
cyclic and separating for $\mathcal{L}(\mathfrak{H})$ (Theorem
{\ref{cyclicth}}), one can construct the modular operator $\Delta_\xi$ (see
\cite{BR}). Namely, if $S$ and $F$ are closures of antilinear operators
given by
\begin{eqnarray*}
S(a\xi)=a^*\xi\textit{ for all }a\in\mathcal{L}(\mathfrak{H})\textit{ and }
F(\xi a^\prime)=\xi\left(a^\prime\right)^*\textit{ for all }a^\prime\in
\mathcal{R}(\mathfrak{H}),
\end{eqnarray*}
then
\begin{eqnarray*}
F=S^*\textit{ and }\;\;\Delta_\xi=FS.
\end{eqnarray*}
Hence, with $a=a_1\otimes a_2\otimes\ldots\otimes a_i\otimes I_2\otimes
I_2\otimes\ldots$ one has
\begin{eqnarray*}
a^*\xi=\xi\cdot\left(\bigotimes\limits_{j=1}^i
\mathfrak{I}_j(\mathbb{I})\right)^{-1}\otimes I_2\otimes
I_2\otimes\ldots\cdot a^*\cdot\left(\bigotimes\limits_{j=1}^i
\mathfrak{I}_j(\mathbb{I})\right)\otimes I_2\otimes I_2 \otimes\ldots.
\end{eqnarray*}
Therefore,
\begin{eqnarray*}
\Delta_\xi\left(a\xi\right)=F\left(a^*\xi\right)=\xi\cdot
\left(\bigotimes\limits_{j=1}^i\mathfrak{I}_j(\mathbb{I})\right)^*\otimes
I_2\otimes\ldots\cdot a\cdot\left(\bigotimes\limits_{j=1}^i
\left(\mathfrak{I}_j(\mathbb{I})\right)^*\right)^{-1}\otimes
I_2\otimes\ldots.
\end{eqnarray*}
Finally, use the relation
$\mathfrak{I}_j(\mathbb{I})\left(\mathfrak{I}_j(\mathbb{I})\right)^*=
\mathfrak{I}_j\mathcal{O}_j^{(0,j)}\mathfrak{I}_j^{-1}$ (see {(\ref{lr})})
to obtain
\begin{eqnarray}\label{modular}
\Delta_\xi\left(a\xi\right)=\bigotimes\limits_{j=1}^i
\left(\mathfrak{I}_j\mathcal{O}_j^{(0,j)}\mathfrak{I}_j^{-1}\right)a
\left(\bigotimes\limits_{j=1}^i\mathfrak{I}_j\mathcal{O}_j^{(0,j)}
\mathfrak{I}_j^{-1}\right)^{-1}\otimes\mathfrak{I}_{i+1}(\mathbb{I})\otimes
\mathfrak{I}_{i+2}(\mathbb{I})\otimes\ldots.
\end{eqnarray}

Thus the modular operator $\Delta_\xi$ is defined in a natural way by the
Okounkov operator $\mathcal{O}_j$ (see (\ref{Cesaro}), \cite{Ok1},
\cite{Ok2}).

\section{The properties of Olshanski's semigroup.}

\begin{Lm}\label{eight}
Let $i(p)$ be an element of $p\in\mathbb{N}\diagup s$. Given any
$\gamma=\left(\gamma_1,\gamma_2,\cdots,\gamma_n,\cdots\right)\in
\Gamma_e^\infty$, there exists $\tilde{\gamma}\in\Gamma_e^\infty$ with the
property $\tilde{\gamma}\cdot
s\cdot\gamma\cdot\tilde{\gamma}^{-1}=s\cdot\gamma^\prime$, where
\begin{eqnarray*}
&\gamma^\prime_{s^{(l-1)}\left(i(p)\right)}=e_{\Gamma}\textit{ for all }
l=1,2,\ldots, |p|-1 \textit{ and }p\in\mathbb{N}\diagup s,
\\ &\gamma^\prime_{s^{(|p|-1)}\left(i(p)\right)}=
\gamma_{s^{(|p|-1)}(i(p))}\cdot\gamma_{s^{(|p|-2)}(i(p))}\cdots
\gamma_{i(p)}.
\end{eqnarray*}
\end{Lm}
\begin{proof}
Let the $\tilde{\gamma}$ be defined as follows:
\begin{eqnarray*}
\tilde{\gamma}_{i(p)}=e_{\Gamma},\tilde{\gamma}_{s(i(p))}=
\gamma_{i(p)}^{-1},\tilde{\gamma}_{s^{(2)}(i(p))}=
\gamma_{i(p)}^{-1}\cdot\gamma_{s(i(p))}^{-1},\ldots,
\\ \tilde{\gamma}_{s^{(|p|-1)}(i(p))}=\gamma_{i(p)}^{-1}\cdot
\gamma_{s(i(p))}^{-1}\cdots\gamma_{s^{(|p|-2)}(i(p))}\textit{ for all
}p\in\mathbb{N}\diagup s.
\end{eqnarray*}
Now our statement can be readily verified.
\end{proof}

\begin{Lm}\label{contigu}
Let $s$ be a cycle from $\mathfrak{S}_\infty$. Suppose that for
$\beta,\gamma\in\Gamma_e^\infty$ the following relations hold:
\begin{eqnarray*}
\beta_k=\gamma_k=e_{\Gamma}\textit{ for all
}k\in\left\{j\in\mathbb{N}\big|\;s(j)=j\right\}.
\end{eqnarray*}
If $s\beta$ and $s\gamma$ are in the same conjugate class, then there
exists $\tilde{\gamma}\in\Gamma_e^\infty$ such that
$s\gamma=\tilde{\gamma}\cdot s\beta\cdot\tilde{\gamma}^{-1}$.
\end{Lm}
\begin{proof}
One may assume without loss of generality that
\begin{eqnarray*}
s(k)=k+1\textit{ for }k=1,2,\ldots,m-1,\;s(m)=1\textit{ and }s(l)=l\textit{
for all }l>m .
\end{eqnarray*}
By Lemma {\ref{eight}} there exist
$\tilde{\gamma},\tilde{\beta}\in\Gamma_e^\infty$ with the properties
\begin{eqnarray}\label{pr}
\begin{split}
&\tilde{\gamma}\cdot
s\cdot\gamma\cdot\tilde{\gamma}^{-1}=s\cdot\gamma^\prime,\;
\tilde{\beta}\cdot s\cdot\beta\cdot\tilde{\beta}^{-1}=s\cdot\beta^\prime,
\textit{ where }
\\ &\gamma^\prime_k=\beta^\prime_k=e_{\Gamma}\textit{ for }k = 1,2,\ldots,
m-1,m+1,\ldots.
\end{split}
\end{eqnarray}
Let $s\in\mathfrak{S}_\infty$ and $\delta\in\Gamma_e^\infty$ be such that
\begin{eqnarray*}
\left(t\delta\right)s\gamma^\prime\left(t\delta\right)^{-1}=s\beta^\prime.
\end{eqnarray*}
One has the following relations:
\begin{eqnarray}\label{relation}
\begin{matrix}
\delta_2\gamma_1^\prime&=&\beta_{t(1)}^\prime\delta_1
\\ \delta_3\gamma_2^\prime&=&\beta_{t(2)}^\prime\delta_2
\\ \cdots&\cdots &\cdots
\\ \delta_m\gamma_{m-1}^\prime&=&\beta_{t(m-1)}^\prime\delta_{m-1}
\\ \delta_1\gamma_{m}^\prime&=&\beta_{t(m)}^\prime\delta_{m}.
\end{matrix}
\end{eqnarray}
By assumptions of the Lemma, $t\left( \left\{ 1,2,\ldots,m \right\}
\right)=\left\{ 1,2,\ldots,m \right\}$, and we may assume that $t(k)=k
\textit{ for all } k>m$. Hence, there exists a map $f$ from $\mathbb{N}$ to
$\mathbb{N}$ such that
\begin{eqnarray*}
t(k)=s^{f(k)} (k) \textit{ for } k\in\mathbb{N}.
\end{eqnarray*}
Now use the relation $ts=st$ to obtain
\begin{eqnarray}
f(k)=l \textit{ for } k=1,2,\ldots,m.
\end{eqnarray}
Since $s^m$ is the identity, it suffices to consider the case
$l\in\left\{1,2,\ldots,m-1 \right\}$.

Use ({\ref{relation}}) to obtain
\begin{eqnarray*}
\delta_1=\ldots=\delta_{m-l},\;\delta_{m-l+1}=\ldots=\delta_m,
\\ \beta_m^\prime=\delta_m\delta_1^{-1},\;\gamma_m^\prime=
\delta_1^{-1}\delta_m.
\end{eqnarray*}
These relations together with ({\ref{pr}}) yield the following relation:
\begin{eqnarray*}
\delta^\prime\;s\gamma^\prime\;\left(\delta^\prime\right)^{-1}=
s\beta^\prime,\textit{ where
}\delta^\prime=\left(\delta_m^{-1}\delta_1,\delta_m^{-1}\delta_1,\ldots,
\delta_m^{-1}\delta_1,\ldots\right).
\end{eqnarray*}
\end{proof}

In what follows, $\left(\pi_\phi,\mathcal{H}_\phi,\xi_\phi\right)$ is the
unitary representation of $G=\Gamma\wr\mathfrak{S}_\infty$ that corresponds
by GNS-construction to the character $\phi$. In particular, the operators
$\pi\left(G\right)$ act in $\mathcal{H}_\phi$ with {\it cyclic} {\it
separating} vector $\xi_\phi$. That is,
\begin{eqnarray}\label{standart}
\left[\pi_\phi\left(G\right)\xi_\phi\right]=\left[\pi_\phi
\left(G\right)^\prime\xi_\phi\right]=\mathcal{H}_\phi,
\end{eqnarray}
where $\left[\mathcal{S}\right]$ stands for the closed subspace in
$\mathcal{H}_\phi$ generated by $\mathcal{S}$. Moreover
$\phi(g)=\left(\pi_\phi\left(g\right)\xi_\phi,\xi_\phi\right)$ for all
$g\in G$.

The property ({\ref{standart}}) allows one to produce a unitary spherical
representation $\pi_\phi^{(2)}$ of the Olshanski pair $\left(G\times
G,K\right)$, where $K=diag\,G=\left\{(g,g)\right\}_{g\in G}$. Namely,
\begin{eqnarray}\label{piphi2}
\pi_\phi^{(2)}\left(g_1,g_2\right)x\xi_\phi=\pi_\phi\left(g_1\right)x
\pi_\phi\left(g_2\right)^*\xi_\phi\;\text{ for all
}x\in\pi_\phi\left(G\right)^{\prime\prime}.
\end{eqnarray}

Let
\begin{eqnarray*}
& G_n(\infty)=\left\{g=s\cdot\gamma\in G\big|\;s(l)=l\textit{ and }
\gamma_l=e\textit{ for all }l=1,2,\cdots,n\right\},
\\ &K_n(\infty)=K\cap\left(G_n(\infty)\times G_n(\infty)\right),G_n=
\Gamma\wr\mathfrak{S}_n,K_n=\left(G_n\times G_n\right)\cap K.
\end{eqnarray*}
It follows from the definition that $G_0(\infty)=G_\infty=G$,
$K_0(\infty)=K_\infty=K$.

Set
\begin{eqnarray*}
\mathcal{H}_\phi^{K_n(\infty)}=\left\{\eta\in\mathcal{H}_\phi\big|\;
\pi_\phi^{(2)}\left(g\right)\eta=\eta\textit{ for all }g\in
K_n(\infty)\right\},
\end{eqnarray*}
and let $P_n$ be an orthogonal projection onto
$\mathcal{H}_\phi^{K_n(\infty)}$.
\begin{Lm}\label{tm}
$\bigcup\limits_{n=0}^\infty\mathcal{H}_\phi^{K_n(\infty)}$ is a dense
subspace in $\mathcal{H}_\phi$. In different terms,
$\lim\limits_{n\to\infty}P_n=\mathcal{I}_{\mathcal{H}_\phi}$ in the strong
operator topology.
\end{Lm}
\begin{proof}
It follows from the definition of $\pi_\phi^{(2)}$ (see ({\ref{piphi2}}))
that
\begin{eqnarray}\label{dense}
\left[\pi_\phi^{(2)}\left(G_n\right)\xi_\phi\right]\subset
\mathcal{H}_\phi^{K_n(\infty)}.
\end{eqnarray}
On the other hand, $\xi_\phi$ is a cyclic vector. That is,
$\left[\bigcup\limits_{n=1}^\infty\pi_\phi^{(2)}
\left(G_n\right)\xi_\phi\right]=\mathcal{H}_\phi$. Now our statement
follows from ({\ref{dense}}).
\end{proof}

Remind a construction of asymptotic operators as it appears in \cite{Ok1},
\cite{Ok2}. Consider the transposition $(i,n)\in\mathfrak{S}_\infty$ and
the operator
\begin{eqnarray}\label{Cesaro}
\mathcal{O}_k=\lim\limits_{n\to\infty}\frac{1}{n}\sum\limits_{l=1}^n
\pi_\phi\left((k,l)\right).
\end{eqnarray}
The limit exists in the strong operator topology.

We follow the idea of Olshanski (see \cite{Ok1}, \cite{O2} and \cite{N1}
for the case of $\mathfrak{S}_\infty$) in extending his approach to our
setting for $\Gamma\wr \mathfrak{S}_\infty$. An algebraic structure of the
associated Olshanski semigroup has been used above to predict important
relations between the operators $\mathcal{O}_k$ and
$\pi_\phi\left(\Gamma_e^\infty\right)$.

Below we sketch the basic algebraic constructions and hope, with some
details being left to a reader. Hopefully, this will allow later to receive
a complete classification of {\it admissible} representations for wreath
products $\Gamma\wr\mathfrak{S}_\infty$ in some reasonable cases.


For any $n\in\mathbb{N}$ consider the $W^*-$algebra
$\left(P_n\pi_\phi^{(2)}\left(G\times G\right)P_n\right)^{\prime\prime}$
generated by the operators $P_n\pi_\phi^{(2)}\left(G\times G\right)P_n$,
which act in $P_n\mathcal{H}_\phi$. Obviously, the map
\begin{eqnarray}\label{reprcoset}
\pi_\phi^{(2,n)}:\;\left(g_1,g_2\right)\in G\times G\to
P_n\pi_\phi^{(2)}\left(g_1,g_2\right)P_n\in\left(P_n\pi_\phi^{(2)}
\left(G\times G\right)P_n\right)^{\prime\prime}
\end{eqnarray}
is constant on the double cosets $K_n(\infty)\diagdown G\times G\diagup
K_n(\infty)$.

We are about to equip $K_n(\infty)\diagdown G\times G\diagup K_n(\infty)$
with a structure of semigroup in such a way that $\pi_\phi^{(2,n)}$ becomes
a homomorphism. First, we extend the idea of Olshanski who applied a
diagram technic to describe algebraic operations on $\mathfrak{S}_\infty$,
onto the semigroup $K_n(\infty)\diagdown G\times G\diagup K_n(\infty)$. For
that, for any double coset we construct a so-called an {\it admissible
graph}, which carries an important information about this coset.

Let $\mathfrak{G}$ be any graph. Denote by $V\left(\mathfrak{G}\right)$ and
$E\left(\mathfrak{G}\right)$ respectively the set of all vertices and edges
of $\mathfrak{G}$. Consider a disjoint union
\begin{eqnarray*}
\overline{\left(\mathbb{Z}\setminus
0\right)}\sqcup\underline{\left(\mathbb{Z}\setminus
0\right)}=\left\{\cdots,\overline{-2},\overline{-1},\overline{1},
\overline{2},\cdots\right\}\sqcup\left\{\cdots,\underline{-2},
\underline{-1},\underline{1},\underline{2},\cdots\right\}
\end{eqnarray*}
of two copies of $\left(\mathbb{Z}\setminus 0\right)$. A map
$\psi_{\mathfrak{G}}$ from $V\left(\mathfrak{G}\right)$ to this disjoint
union such that $\psi_{\mathfrak{G}}(v)\ne\psi_{\mathfrak{G}}(v^\prime)$
for all pairs of different vertices $v$ and $v^\prime$, is called a {\it
vertex-coloring} of the graph $\mathfrak{G}$.
\begin{Def}\label{adgraph}
An {\it admissible graph} is a {\it vertex-colored}, {\it directed} graph
$\mathfrak{G}$ with the properties:
\begin{itemize}
\item [\it (a)] $V\left(\mathfrak{G}\right)$ is the disjunct union of four
finite sets $V\left(\mathfrak{G}\right)_u^-$,
$V\left(\mathfrak{G}\right)_u^+$, $V\left(\mathfrak{G}\right)_o^-$ and
$V\left(\mathfrak{G}\right)_o^+$, with
\begin{eqnarray*}
\psi_{\mathfrak{G}}\left(V\left(\mathfrak{G}\right)_u^\pm\right)=
\left\{\underline{\pm 1},\underline{\pm 2},\cdots,\underline{\pm
n}\right\},\psi_{\mathfrak{G}}\left(V\left(\mathfrak{G}\right)_o^\pm\right)=
\left\{\overline{\pm 1},\overline{\pm 2},\cdots,\overline{\pm n}\right\};
\end{eqnarray*}
\item [\it (b)] if $i(e)$ and $t(e)$ are {\it initial} and {\it terminal}
vertices of $e\in E\left(\mathfrak{G}\right)$, then
\begin{eqnarray*}
\left|\left\{ e\in E\left(\mathfrak{G}\right)\Big|\;i(e)=v,t(e)=v^\prime
\textit{ or }t(e)=v,i(e)=v^\prime\right\}\right|=\left\{
\begin{array}{rl}
1&\textit{ if }v\ne v^\prime
\\ 0&\textit{ otherwise};
\end{array}\right.
\end{eqnarray*}
\item [\it (c)] if $e\in E\left(\mathfrak{G}\right)$ then either of the
four cases holds:
\begin{eqnarray*}
i(e)\in V\left(\mathfrak{G}\right)_o^-\textit{ and }\;t(e)\in
V\left(\mathfrak{G}\right)_u^-,
\\ i(e)\in V\left(\mathfrak{G}\right)_u^+\textit{ and }\;t(e)\in
V\left(\mathfrak{G}\right)_u^-,
\\ i(e)\in V\left(\mathfrak{G}\right)_u^+\textit{ and }\;t(e)\in
V\left(\mathfrak{G}\right)_o^+,
\\ i(e)\in V\left(\mathfrak{G}\right)_o^-\textit{ and }\;t(e)\in
V\left(\mathfrak{G}\right)_o^+;
\end{eqnarray*}
\item [\it (d)] one has a well defined {\it marking} function
$m_\mathfrak{G}$ from $E\left(\mathfrak{G}\right)$ to
$\left\{\mathbb{N}\diagup 2,\;0\right\}\times\Gamma$, where
$\mathbb{N}\diagup 2$ is the set of positive half-integer numbers.
\end{itemize}
\end{Def}
\begin{Def}\label{diag}
The disjunct sum of the {\it admissible} graph $\mathfrak{G}$ and the
countable set $\mathfrak{C}$ of the circles is called the {\it admissible
diagram} if there is a well defined {\it marking} function $m_\mathfrak{C}$
from $\mathfrak{C}$ to $\mathbb{N}\times\mathfrak{c}_{_\Gamma}(\Gamma)$,
where $\mathfrak{c}_{_\Gamma}(\Gamma)$ is the set of conjugacy classes of
$\Gamma$ with property: there are finitely many elements
$m_\mathfrak{C}(c)\;\left(c\in\mathfrak{C}\right)$ not equal to
$\left(1,e_\Gamma\right)$.
\end{Def}

Here is a graphic interpretation for the elements of $G\times G$. Let
$g=\left(s_1\gamma^\prime,s_2\gamma^{\prime\prime}\right)$. Consider two
copies $ \overline{\left(\mathbb{Z}\setminus 0\right)}$ and
$\underline{\left(\mathbb{Z}\setminus 0\right)}$ of $\mathbb{Z}\setminus
0$. It is convenient to position the elements of
$\overline{\left(\mathbb{Z}\setminus 0\right)}$ and
$\underline{\left(\mathbb{Z}\setminus 0\right)}$ on two horizontal lines,
$\overline{\left(\mathbb{Z}\setminus 0\right)}$ above $
\underline{\left(\mathbb{Z}\setminus 0\right)}$. Draw the edges
$e\left(\underline{i},\overline{s_2(i)}\right)$ from
$\underline{i}\in\underline{\left(\mathbb{Z}\setminus 0\right)}$ to
$\overline{s_2(i)}\in\overline{\left(\mathbb{Z}\setminus 0\right)}$ and
$e\left(\underline{-i},\overline{-s_1(i)}\right)$ from
$\underline{-i}\in\underline{\left(\mathbb{Z}\setminus 0\right)}$ to
$\overline{-s_1(i)}\in\overline{\left(\mathbb{Z}\setminus 0\right)}$ for
$i>0$ (see Fig. \ref{pic1}). Finally, define in a similarity with
Definition {\ref{adgraph}} the {\it marking} function
$m_{\mathfrak{gr}(g)}$ as follows:
\begin{eqnarray*}\label{mark}
m_{\mathfrak{gr}(g)}\left(e\left(\underline{i},\overline{j}\right)\right)=
\left\{
\begin{array}{rl}
\left(0,\gamma_i^{\prime\prime}\right)&\textit{ if }i>0,
\\ \left(0,\gamma_{|i|}^{\prime}\right)&\textit{ if }i<0.
\end{array}\right.
\end{eqnarray*}

\begin{figure}
\caption{$\mathfrak{gr}(g)$}\label{pic1}
 \begin{picture}(300,85)(0,30)
\multiput(300,100)(-30,0){10}{$\bullet$}
\multiput(300,40)(-30,0){10}{$\bullet$}
 \put(300,107){$\overline{5}$}
\put(300,33){$\underline{5}$}
 \put(270,107){$\overline{4}$}
\put(270,33){$\underline{4}$}
 \put(240,107){$\overline{3}$}
\put(240,33){$\underline{3}$}
 \put(210,107){$\overline{2}$}
\put(210,33){$\underline{2}$}
 \put(180,107){$\overline{1}$}
\put(180,33){$\underline{1}$}
 \put(150,107){$\overline{-1}$}
\put(150,33){$\underline{-1}$}
 \put(120,107){$\overline{-2}$}
\put(120,33){$\underline{-2}$}
 \put(90,107){$\overline{-3}$}
\put(90,33){$\underline{-3}$}
 \put(60,107){$\overline{-4}$}
\put(60,33){$\underline{-4}$}
 \put(30,107){$\overline{-5}$}
\put(30,33){$\underline{-5}$}

 \put(31,42){\vector(2,1){120}}
 \put(61,42){\vector(1,1){60}}
 \put(90,44){\vector(-1,2){28}}
 \put(120,43){\vector(-3,2){86}}
 \put(150,43){\vector(-1,1){57}}
 \put(182,42){\vector(3,2){89}}
 \put(212,42){\vector(3,2){89}}
\put(240,44){\vector(-1,2){28}}
 \put(270,44){\vector(-1,2){28}}
\put(300,43){\vector(-2,1){115}}

 \put(285,53){$\gamma_5^{\prime\prime}$}
 \put(253,50){$\gamma_4^{\prime\prime}$}
  \put(217,90){$\gamma_3^{\prime\prime}$}
  \put(205,66){$\gamma_1^{\prime\prime}$}
  \put(275,80){$\gamma_2^{\prime\prime}$}
  \put(128,67){$\gamma_1^{\prime}$}
 \put(105,54){$\gamma_2^{\prime}$}
 \put(70,85){$\gamma_3^{\prime}$}
 \put(70,46){$\gamma_4^{\prime}$}
 \put(46,57){$\gamma_5^{\prime}$}
\end{picture}
\end{figure}

A graph has been constructed, to be denoted by $\mathfrak{gr}(g)$.
Obviously, $i\left(E\left(\mathfrak{gr}(g)\right)\right)=
\underline{\left(\mathbb{Z}\setminus 0\right)}$ and
$t\left(E\left(\mathfrak{gr}(g)\right)\right)=
\overline{\left(\mathbb{Z}\setminus 0\right)}$.

To produce the graph $\mathfrak{gr}(gh)$, it is convenient to position
$\mathfrak{gr}(g)$ above $\mathfrak{gr}(h)$. After the natural gluing the
vertices $\overline{i}\in t\left(E\left(\mathfrak{gr}(h)\right)\right)$ and
$\underline{i}\in i\left(E\left(\mathfrak{gr}(g)\right)\right)$ we receive
$\mathfrak{gr}(gh)$. It is clear that
\begin{eqnarray*}
m_{\mathfrak{gr}(gh)}\left(e\left(\underline{i},\overline{j}\right)\right)=
m_{\mathfrak{gr}(g)}\left(e_g\left(\underline{k},\overline{j}\right)\right)
\cdot m_{\mathfrak{gr}(h)}\left(e_h\left(\underline{i},\overline{k}\right)
\right),
\end{eqnarray*}
where the edge $e\left(\underline{i},\overline{j}\right)$ is the joining of
the edges $e_h\left(\underline{i},\overline{k}\right)\in
E\left(\mathfrak{gr}(h)\right)$ and
$e_g\left(\underline{k},\overline{j}\right)\in
E\left(\mathfrak{gr}(g)\right)$. If $h\in G\times G$ is defined by the
graph on Fig. \ref{pic2}, then the graph $\mathfrak{gr}(gh)$ on Fig.
\ref{pic3} corresponds to the product $gh$. Besides that, we equip
$\left\{\mathbb{N}\diagup 2,\;0\right\}\times\Gamma$ with a natural
semigroup structure.

\begin{figure}
\caption{$\mathfrak{gr}(h)$}\label{pic2}
\begin{picture}(300,85)(0,30)
\multiput(300,100)(-30,0){10}{$\bullet$}
\multiput(300,40)(-30,0){10}{$\bullet$}
 \put(300,107){$\overline{5}$}
\put(300,33){$\underline{5}$}
 \put(270,107){$\overline{4}$}
\put(270,33){$\underline{4}$}
 \put(240,107){$\overline{3}$}
\put(240,33){$\underline{3}$}
 \put(210,107){$\overline{2}$}
\put(210,33){$\underline{2}$}
 \put(180,107){$\overline{1}$}
\put(180,33){$\underline{1}$}
 \put(150,107){$\overline{-1}$}
\put(150,33){$\underline{-1}$}
 \put(120,107){$\overline{-2}$}
\put(120,33){$\underline{-2}$}
 \put(90,107){$\overline{-3}$}
\put(90,33){$\underline{-3}$}
 \put(60,107){$\overline{-4}$}
\put(60,33){$\underline{-4}$}
 \put(30,107){$\overline{-5}$}
\put(30,33){$\underline{-5}$}

 \put(32,42){\vector(3,2){89}}
 \put(63,42){\vector(1,2){29}}
 \put(93,44){\vector(1,1){58}}
 \put(120,43){\vector(-1,1){58}}
 \put(150,43){\vector(-2,1){116}}

 \put(182,42){\vector(2,1){119}}
 \put(212,42){\vector(1,1){59}}
\put(240,44){\vector(-1,1){58}}
 \put(272,44){\vector(-1,2){29}}
\put(300,43){\vector(-3,2){87}}

 \put(285,54){$\delta_5^{\prime\prime}$}
 \put(256,50){$\delta_4^{\prime\prime}$}
  \put(198,87){$\delta_3^{\prime\prime}$}
  \put(253,95){$\delta_2^{\prime\prime}$}
  \put(275,82){$\delta_1^{\prime\prime}$}

  \put(135,79){$\delta_3^{\prime}$}
 \put(50,82){$\delta_1^{\prime}$}
 \put(72,92){$\delta_2^{\prime}$}
 \put(70,50){$\delta_4^{\prime}$}
 \put(46,60){$\delta_5^{\prime}$}
\end{picture}
\end{figure}

\begin{figure}
\caption{$\mathfrak{gr}(gh)$}\label{pic3}

 \begin{picture}(305,85)(0,30)
\multiput(300,100)(-30,0){10}{$\bullet$}
\multiput(300,40)(-30,0){10}{$\bullet$}
 \put(300,107){$\overline{5}$}
\put(300,33){$\underline{5}$}
 \put(270,107){$\overline{4}$}
\put(270,33){$\underline{4}$}
 \put(240,107){$\overline{3}$}
\put(240,33){$\underline{3}$}
 \put(210,107){$\overline{2}$}
\put(210,33){$\underline{2}$}
 \put(180,107){$\overline{1}$}
\put(180,33){$\underline{1}$}
 \put(150,107){$\overline{-1}$}
\put(150,33){$\underline{-1}$}
 \put(120,107){$\overline{-2}$}
\put(120,33){$\underline{-2}$}
 \put(90,107){$\overline{-3}$}
\put(90,33){$\underline{-3}$}
 \put(60,107){$\overline{-4}$}
\put(60,33){$\underline{-4}$}
 \put(30,107){$\overline{-5}$}
\put(30,33){$\underline{-5}$}

 \put(32,42){\vector(0,1){58}}
 \put(62,42){\vector(0,1){58}}
 \put(92,42){\vector(0,1){58}}
 \put(122,43){\vector(0,1){58}}
 \put(152,43){\vector(0,1){58}}

 \put(182,43){\vector(0,1){58}}
 \put(212,42){\vector(1,2){29}}
\put(242,42){\vector(1,2){29}}
 \put(272,44){\vector(-1,1){58}}
\put(303,43){\vector(0,1){58}}

 \put(303,70){$\gamma_2^{\prime\prime}\delta_5^{\prime\prime}$}
 \put(263,54){$\gamma_3^{\prime\prime}\delta_4^{\prime\prime}$}
  \put(183,83){$\gamma_5^{\prime\prime}\delta_3^{\prime\prime}$}

  \put(200,60){$\gamma_4^{\prime\prime}
\delta_2^{\prime\prime}$}

  \put(265,82){$\gamma_1^{\prime\prime}
\delta_3^{\prime\prime}$}

  \put(133,79){$\gamma_5^{\prime}\delta_1^{\prime}$}
 \put(103,70){$\gamma_4^{\prime}\delta_2^{\prime}$}
 \put(74,70){$\gamma_1^{\prime}\delta_3^{\prime}$}
 \put(44,70){$\gamma_3^{\prime}\delta_4^{\prime}$}
 \put(14,60){$\gamma_2^{\prime}\delta_5^{\prime}$}
\end{picture}
\end{figure}

Pass to a construction of the {\it admissible diagram} (see Definition
{\ref{diag}}) which corresponds to the coset $\theta_n(g)\in
K_n(\infty)\diagdown G\times G\diagup K_n(\infty)$ containing $g\in G\times
G$. It splits into four steps.
\begin{itemize}
\item At Step 1, draw in $\mathfrak{gr}(g)$ for $i>n$ the edges
$e\left(\underline{-i},\underline{i}\right)$ and
$e\left(\overline{i},\overline{-i}\right)$ that connect the vertices
$\underline{-i}$ to $\underline{i}$ and $\overline{i}$ to $\overline{-i}$.
Denote the new graph by $\overline{\mathfrak{gr}(g)}$ (see Fig.
{\ref{pic4}}) and extend the {\it marking} function $m_{\mathfrak{gr}(g)}$
on $\overline{\mathfrak{gr}(g)}$ assuming that
\begin{eqnarray*}
m_{\overline{\mathfrak{gr}(g)}}
\left(e\left(\underline{-i},\underline{i}\right)\right)=
m_{\overline{\mathfrak{gr}(g)}}
\left(e\left(\overline{i},\overline{-i}\right)\right)=
\left(\frac{1}{2},e_\Gamma\right)\in\left\{\mathbb{N}\diagup
2,\;0\right\}\times\Gamma.
\end{eqnarray*}
\begin{figure}
\caption{$\overline{\mathfrak{gr}(g)}$}\label{pic4}
 \begin{picture}(300,145)
\multiput(300,100)(-30,0){10}{$\bullet$}
\multiput(300,40)(-30,0){10}{$\bullet$}

 \put(240,106){$\overline{3}$}
\put(240,33){$\underline{3}$}
 \put(210,107){$\overline{2}$}
\put(210,33){$\underline{2}$}
 \put(180,107){$\overline{1}$}
\put(180,33){$\underline{1}$}
 \put(150,107){$\overline{-1}$}
\put(150,33){$\underline{-1}$}
 \put(120,107){$\overline{-2}$}
\put(120,33){$\underline{-2}$}
 \put(90,107){$\overline{-3}$}
\put(90,33){$\underline{-3}$}

 \put(150,101){\vector(-2,-1){115}}
 \put(121,101){\vector(-1,-1){57}}
 \put(33,101){\vector(3,-2){88}}
 \put(63,101){\vector(1,-2){29}}
 \put(152,42){\vector(-1,1){59}}
 \put(182,42){\vector(3,2){89}}
 \put(212,42){\vector(3,2){89}}
\put(240,44){\vector(-1,2){28}}
 \put(270,44){\vector(-1,2){28}}
\put(300,43){\vector(-2,1){115}}

\spline(300,40)(270,20)(165,0)(60,20)(32,42)
\spline(270,40)(240,20)(165,5)(90,20)(60,42)
\put(220,5){$\scriptscriptstyle{\left( \frac{1}{2},0\right)}$}
\put(165,15){$\scriptstyle{\left( \frac{1}{2},0\right)}$}

\spline(305,102)(270,120)(165,140)(60,120)(32,102)
\spline(275,102)(240,120)(165,135)(90,120)(62,102)
\put(220,133){$\scriptstyle{\left( \frac{1}{2},0\right)}$}
\put(165,123){$\scriptstyle{\left( \frac{1}{2},0\right)}$}

\put(285,53){$\scriptstyle{(0,\gamma_5^{\prime\prime})}$}
 \put(243,50){$\scriptstyle{(0,\gamma_4^{\prime\prime})}$}
  \put(219,90){$\scriptstyle{(0,\gamma_3^{\prime\prime})}$}
  \put(193,66){$\scriptstyle{(0,\gamma_1^{\prime\prime})}$}
  \put(275,80){$\scriptstyle{(0,\gamma_2^{\prime\prime})}$}
  \put(128,67){$\scriptstyle{(0,\gamma_1^{\prime})}$}

\dottedline{1}(50,90)(55,123)(48,130)
 \put(10,132){\fbox{$\scriptstyle{\left(0,
 \left(\gamma_2^\prime\right)^{-1}\right)}$}}


 \dottedline{1}(75,55)(80,40)(75,25)(50,10)(37,10)
\put(0,10){\fbox{$\scriptstyle{\left(0,
 \left(\gamma_4^\prime\right)^{-1}\right)}$}}

 \put(20,56){$\scriptstyle{(0,(\gamma_5^{\prime})^{-1})}$}

 \dottedline{1}(70,85)(38,80)
 \put(0,75){\fbox{$\scriptstyle{\left(0,
 \left(\gamma_3^\prime\right)^{-1}\right)}$}}
\end{picture}

\end{figure}
\item At Step 2, extend the {\it marking} function $m_{\mathfrak{gr}(g)}$
to the pathes of $\overline{\mathfrak{gr}(g)}$. First, it is reasonable to
assume, that
\begin{eqnarray*}
m_{\overline{\mathfrak{gr}(g)}}
\left(e\left(\overline{j},\underline{i}\right)\right)=
\left(l,\gamma^{-1}\right),\textit{ when }m_{\overline{\mathfrak{gr}(g)}}
\left(e\left(\underline{i},\overline{j}\right)\right)=(l,\gamma).
\end{eqnarray*}
If the path $\mathfrak{p}\in\overline{\mathfrak{gr}(g)}$ is formed by a
sequence $\left\{e_1,e_2,\ldots,e_j\right\}$ of edges, then
\begin{eqnarray*}
m_{\overline{\mathfrak{gr}(g)}}(\mathfrak{p})=
m_{\overline{\mathfrak{gr}(g)}}\left(e_j\right)\cdot
m_{\overline{\mathfrak{gr}(g)}}\left(e_{j-1}\right)\cdots
m_{\overline{\mathfrak{gr}(g)}}\left(e_1\right).
\end{eqnarray*}
\item At Step 3, notice, that $\overline{\mathfrak{gr}(g)}$ is a disjoint
union of its connected components, each of those being a {\it non closed}
path or {\it cycle}. By the construction, the ends of any non-closed path
belong to
\begin{eqnarray*}
\mathfrak{V}_n=\left\{\overline{-n},\cdots,\overline{-2},\overline{-1},
\overline{1},\overline{2},\cdots,\overline{n}\right\}\sqcup
\left\{\underline{-n},\cdots,\underline{-2},\underline{-1},\underline{1},
\underline{2},\cdots,\underline{n}\right\}.
\end{eqnarray*}
Furthermore, define a coherent {\it positive} orientation on any non-closed
path, which contains vertices $\overline{i}$ or $\underline{i}$ with
$|i|>n$, if we assume that its initial vertex belongs to the set
\begin{eqnarray*}
\left\{\cdots,\overline{-2},\overline{-1}\right\}\sqcup
\left\{\underline{1},\underline{2},\cdots \right\}.
\end{eqnarray*}
The condition that any cycle contains an edge of the form
$e\left(\underline{i},\overline{j}\right)$, where $i,j>0$, defines the
corresponding orientation on closed pathes.

\item Step 4. The set of vertices of the {\it admissible} graph
$\mathfrak{gr}\left(\theta_n(g)\right)$ (see Definition {\ref{adgraph}}) is
$\mathfrak{V}_n$. Each oriented non-closed path $\mathfrak{p}$ (see Step 3)
determines an oriented edge $e(\mathfrak{p})$ of
$\mathfrak{gr}\left(\theta_n(g)\right)$ as follows:
\begin{eqnarray*}
i(\mathfrak{p})=i(e(\mathfrak{p})),\;t(\mathfrak{p})=t(e(\mathfrak{p})).
\end{eqnarray*}
Finally, $m_{\mathfrak{gr}\left(\theta_n(g)\right)}(e(\mathfrak{p}))=
m_{\overline{\mathfrak{gr}(g)}}(\mathfrak{p})$ (see Definition
{\ref{adgraph}} ({\it d})). Similarly, each oriented cycle $\mathfrak{c}$
defines the circle $c(\mathfrak{c})$ (see Definition {\ref{diag}}). The
value of the marking function on $c(\mathfrak{c})$ coincides with $
m_{\overline{\mathfrak{gr}(g)}}(\mathfrak{c})=(r,\gamma)$. It is clear that
the conjugacy class of $\gamma$ does not depend on a choice of initial
vertex in $\mathfrak{c}$. The visualization of this algorithm for $g\in
G_5$ is depicted on Fig. {\ref{pic5}}. Notice that there are no circles on
Fig. {\ref{pic5}} which are marked through $\left(1,e_\Gamma\right)$.
\end{itemize}
Hence, by Definition {\ref{diag}}, we obtain the {\it admissible diagram}
$\mathfrak{d}\left(\theta_n(g)\right)$ of the coset $\theta_n(g)$.

Using the diagram technic, we describe the algorithm of multiplication for
cosets. Let $\theta_n\left(g\right)$ and $\theta_n\left(h\right)$ be two
elements from $K_n(\infty)\diagdown G\times G\diagup K_n(\infty)$. Let
$\theta_n\left(g\right)\circ\theta_n\left(h\right)$ stand for the product
$\theta_n\left(g\right)$ and $\theta_n\left(h\right)$. Again, it is
convenient to position $\mathfrak{d}\left(\theta_n\left(g\right)\right)$
above $\mathfrak{d}\left(\theta_n\left(h\right)\right)$. Later on, we paste
the vertices
\begin{eqnarray*}
\left\{\overline{-n},\cdots,\overline{-2},\overline{-1},\overline{1},
\overline{2},\cdots,\overline{n}\right\}\subset\mathfrak{d}
\left(\theta_n\left(h\right)\right)
\end{eqnarray*}
with the corresponding ones from
\begin{eqnarray*}
\left\{\underline{-n},\cdots,\underline{-2},\underline{-1},\underline{1},
\underline{2},\cdots,\underline{n}\right\}\subset\mathfrak{d}
\left(\theta_n\left(g\right)\right).
\end{eqnarray*}
The resulting graph inherits naturally an orientation from the diagrams
$\mathfrak{d}\left(\theta_n\left(g\right)\right)$ and
$\mathfrak{d}\left(\theta_n\left(h\right)\right)$. Just as at Step 4, we
replace the newly formed connected components (non closed pathes or cycles)
by the oriented edges or circles and define on those the marking function.
The received diagram for $g,h\in G_5$ is represented on Fig. {\ref{pic7}}.

\begin{figure}

\caption{${\mathfrak{d}\left(\theta_n(g)\right)}$}\label{pic5}


 \begin{picture}(310,85)(0,30)
\multiput(300,100)(-60,0){6}{$\bullet$}
\multiput(300,40)(-60,0){6}{$\bullet$}
 \put(300,107){$\overline{3}$}

 \spline(122,102)(150,90)(182,102)
 \put(150,92){$\scriptstyle{\blacktriangleright}$}

 \spline(62,102)(180,80)(300,102)
 \put(180,84){$\scriptstyle{\blacktriangleright}$}

 \spline(180,42)(90,60)(0,42)
 \put(60,53){$\scriptstyle{\blacktriangleleft}$}
 \put(122,42){\vector(-2,1){118}}
 \put(303,42){\vector(-1,1){59}}

 \spline(240,42)(150,60)(60,42)
 \put(150,54){$\scriptstyle{\blacktriangleleft}$}

\put(300,33){$\underline{3}$}
 \put(240,107){$\overline{2}$}

\put(240,33){$\underline{2}$}
 \put(185,107){$\overline{1}$}
\put(180,33){$\underline{1}$}

 \put(110,107){$\overline{-1}$}
\put(115,33){$\underline{-1}$}
 \put(60,107){$\overline{-2}$}
\put(60,33){$\underline{-2}$}
 \put(1,107){$\overline{-3}$}
\put(1,33){$\underline{-3}$}
 \put(23,74){$\scriptstyle{\left( 0,
 \gamma_1^\prime\right)}$}

 \put(243,74){$\scriptstyle{\left( 0,
 \gamma_3^{\prime\prime}\right)}$}

 \put(127,105){$\scriptstyle{\left( \frac{1}{2},
 \gamma_5^{\prime\prime}\left(
 \gamma_5^{\prime}\right)^{-1}\right)}$}

\put(92,80){$\scriptstyle{\left( \frac{1}{2},
 \gamma_4^{\prime\prime}\left(
 \gamma_4^{\prime}\right)^{-1}\right)}$}

\put(2,57){$\scriptstyle{\left( \frac{1}{2},
 \left(\gamma_3^\prime\right)^{-1}\gamma_1^{\prime\prime}
\right) }$}

\put(170,60){$\scriptstyle{\left( \frac{1}{2},
 \left(\gamma_2^\prime\right)^{-1}\gamma_2^{\prime\prime}
\right) }$}

\end{picture}

\end{figure}

\begin{figure}
\caption{${\mathfrak{d}\left(\theta_n(h)\right)}$}\label{pic6}

 \begin{picture}(310,85)(0,30)
\multiput(300,100)(-60,0){6}{$\bullet$}
\multiput(300,40)(-60,0){6}{$\bullet$}
 \put(300,107){$\overline{3}$}

 \spline(62,102)(150,90)(242,102)
 \put(150,91){$\scriptstyle{\blacktriangleright}$}

 \spline(2,102)(180,80)(300,102)
 \put(150,85){$\scriptstyle{\blacktriangleright}$}

 \spline(180,42)(150,50)(120,42)
 \put(150,47){$\scriptstyle{\blacktriangleleft}$}
 \put(1,42){\vector(2,1){120}}
 \put(303,42){\vector(-2,1){119}}

 \spline(240,42)(150,60)(60,42)
 \put(150,55){$\scriptstyle{\blacktriangleleft}$}

\put(300,33){$\underline{3}$}
 \put(240,107){$\overline{2}$}

\put(240,33){$\underline{2}$}
 \put(184,107){$\overline{1}$}
\put(184,33){$\underline{1}$}

 \put(112,107){$\overline{-1}$}
\put(110,33){$\underline{-1}$}
 \put(60,107){$\overline{-2}$}
\put(60,33){$\underline{-2}$}
 \put(1,107){$\overline{-3}$}
\put(1,33){$\underline{-3}$}

 \put(128,35){$\scriptstyle{\left( \frac{1}{2},
 \left(\delta_1^\prime\right)^{-1}\delta_1^{\prime\prime}
\right) }$}

 \put(80,60){$\scriptstyle{\left( \frac{1}{2},
 \left(\delta_2^\prime\right)^{-1}\delta_2^{\prime\prime}
\right) }$}

 \put(160,78){$\scriptstyle{\left( \frac{1}{2},
 \delta_4^{\prime\prime}\left(\delta_4^\prime\right)^{-1}
\right) }$}

 \put(128,99){$\scriptstyle{\left( \frac{1}{2},
 \delta_5^{\prime\prime}\left(\delta_5^\prime\right)^{-1}
\right) }$}

 \put(253,70){$\scriptstyle{\left( 0,
 \delta_3^{\prime\prime}
\right) }$}

 \put(30,70){$\scriptstyle{\left( 0,
 \delta_3^\prime
\right) }$}

\end{picture}

\end{figure}

\begin{figure}
\caption{$\mathfrak{d}\left(\theta_n(g)\circ\theta_n(h)\right)$}\label{pic7}
\begin{picture}(310,115)
\multiput(300,100)(-60,0){6}{$\bullet$}
\multiput(300,40)(-60,0){6}{$\bullet$}
 \put(300,107){$\overline{3}$}

 \spline(180,42)(150,50)(120,42)
 \put(150,47){$\scriptstyle{\blacktriangleleft}$}
 \put(2,42){\vector(0,1){60}}
 \put(303,42){\line(-1,1){19}}
  \put(264,80){\vector(-1,1){20}}

 \spline(240,42)(150,60)(60,42)
 \put(150,55){$\scriptstyle{\blacktriangleleft}$}

 \spline(122,102)(150,90)(182,102)
 \put(150,92){$\scriptstyle{\blacktriangleright}$}

 \spline(62,102)(180,80)(300,102)
 \put(180,84){$\scriptstyle{\blacktriangleright}$}

\put(300,33){$\underline{3}$}
 \put(240,107){$\overline{2}$}

\put(240,33){$\underline{2}$}
 \put(184,107){$\overline{1}$}
\put(184,33){$\underline{1}$}

 \put(112,107){$\overline{-1}$}
\put(110,33){$\underline{-1}$}
 \put(60,107){$\overline{-2}$}
\put(60,33){$\underline{-2}$}
 \put(1,107){$\overline{-3}$}
\put(1,33){$\underline{-3}$}

 \put(128,35){$\scriptstyle{\left( \frac{1}{2},
 \left(\delta_1^\prime\right)^{-1}\delta_1^{\prime\prime}
\right) }$}

 \put(80,60){$\scriptstyle{\left( \frac{1}{2},
 \left(\delta_2^\prime\right)^{-1}\delta_2^{\prime\prime}
\right) }$}

\put(3,70){$\scriptstyle{\left( 0,
 \gamma_1^\prime\delta_3^\prime\right)}$}

 \put(127,105){$\scriptstyle{\left( \frac{1}{2},
 \gamma_5^{\prime\prime}\left(
 \gamma_5^{\prime}\right)^{-1}\right)}$}

\put(154,76){$\scriptstyle{\left( \frac{1}{2},
 \gamma_4^{\prime\prime}\left(
 \gamma_4^{\prime}\right)^{-1}\right)}$}

 \put(210,70){\fbox{$\scriptstyle{\left(1,
 \gamma_3^{\prime\prime}\delta_4^{\prime\prime}
 \left(\delta_4^\prime\right)^{-1}
\left( \gamma_3^\prime
 \right)^{-1} \gamma_1^{\prime\prime}
   \delta_3^{\prime\prime}
 \right)}$}}

 \put(30,15){\circle{20}}
 \dottedline{1}(45,14)(70,10)
 \put(43,13){$\scriptscriptstyle{<}$}

 \put(70,13){\fbox{$\scriptstyle{
 \left(1,
\left( \gamma_2^\prime\right)^{-1}
 \gamma_2^{\prime\prime}
 \delta_5^{\prime\prime}  \left(\delta_5^\prime  \right)^{-1}
\right)}$}}
\end{picture}
\end{figure}

\begin{Rem}
Let $ (i,n+1)$ be the transposition from $\mathfrak{S}_\infty\subset
G=\mathfrak{S}_\infty\wr\Gamma$, and let $\gamma$ be any element from
$\Gamma\subset G$. It is easy to check up, following Steps 1-4, that they
are determined by the {\it admissible diagrams} depicted in Figure
{\ref{pic8}}.
\begin{figure}
\caption{$\mathfrak{d}\left(\theta_n((i,n+1))\right)$ and
$\mathfrak{d}\left(\theta_n(\gamma)\right)$}\label{pic8}
 \begin{picture}(320,70)
 \multiput(140,50)(-35,0){5}{$\bullet$}
 \multiput(140,10)(-35,0){5}{$\bullet$}

 \put(140,57){$\overline{i+2}$}
 \put(140,3){$\underline{i+2}$}
 \put(143,11){\vector(0,1){40}}

 \put(105,57){$\overline{i+1}$}
 \put(105,3){$\underline{i+1}$}
 \put(108,11){\vector(0,1){40}}
\put(145,25){$\scriptstyle{\left( 0,e_\Gamma\right)}$}
 \put(110,25){$\scriptstyle{\left( 0,e_\Gamma\right)}$}
 \put(75,25){$\scriptstyle{\left( 1,e_\Gamma\right)}$}
\put(40,25){$\scriptstyle{\left( 0,e_\Gamma\right)}$}
 \put(5,25){$\scriptstyle{\left( 0,e_\Gamma\right)}$}
 \put(70,57){$\overline{i}$}
 \put(70,3){$\underline{i}$}

\put(73,11){\vector(0,1){40}}
 \put(35,57){$\overline{i-1}$}
 \put(35,3){$\underline{i-1}$}
 \put(38,11){\vector(0,1){40}}

 \put(0,57){$\overline{i-2}$}
 \put(0,3){$\underline{i-2}$}
\put(3,11){\vector(0,1){40}}

\multiput(315,50)(-35,0){5}{$\bullet$}
 \multiput(315,10)(-35,0){5}{$\bullet$}

 \put(315,57){$\overline{i+2}$}
 \put(315,3){$\underline{i+2}$}
 \put(318,11){\vector(0,1){40}}

 \put(280,57){$\overline{i+1}$}
 \put(280,3){$\underline{i+1}$}
 \put(283,11){\vector(0,1){40}}
\put(320,25){$\scriptstyle{\left( 0,\gamma_{i+2}\right)}$}
 \put(285,25){$\scriptstyle{\left( 0,\gamma_{i+1}\right)}$}
 \put(250,25){$\scriptstyle{\left( 0,\gamma_i\right)}$}
\put(215,25){$\scriptstyle{\left( 0,\gamma_{i-1}\right)}$}
 \put(180,25){$\scriptstyle{\left( 0,\gamma_{i-2}\right)}$}
 \put(245,57){$\overline{i}$}
 \put(245,3){$\underline{i}$}

\put(248,11){\vector(0,1){40}}
 \put(210,57){$\overline{i-1}$}
 \put(210,3){$\underline{i-1}$}
 \put(213,11){\vector(0,1){40}}

 \put(175,57){$\overline{i-2}$}
 \put(175,3){$\underline{i-2}$}
\put(178,11){\vector(0,1){40}}

 \dottedline{2}(170,0)(170,70)
 \end{picture}

\end{figure}
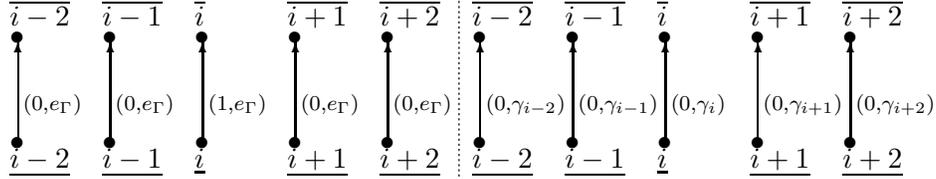
By the graphic interpretation of the multiplication for cosets one has
\begin{eqnarray}
\theta_n(\gamma)\circ\theta_n((i,n+1))=\theta_n((i,n+1))\circ
\theta_n(\gamma).
\end{eqnarray}
If $\pi_\phi^{(2,n)}$ (see (\ref{reprcoset})) is the
representation\label{com} of the semigroup of cosets, then
\begin{eqnarray*}
\pi_\phi^{(2,n)}\left(\theta_n((i,n+1))\right)=P_n\mathcal{O}_iP_n
\text{(see (\ref{Cesaro})) }\;\text{ and
}\pi_\phi^{(2,n)}\left(\theta_n(\gamma)\right)=\pi_\phi^{(2,n)}
\left(\gamma\right).
\end{eqnarray*}
Therefore, $P_n\mathcal{O}_iP_n\pi_\phi^{(2)}(\gamma)P_n=
P_n\pi_\phi^{(2)}(\gamma)P_n\mathcal{O}_iP_n$. Using Lemma {\ref{tm}}, we
obtain
$\mathcal{O}_i\pi_\phi^{(2)}(\gamma)=\pi_\phi^{(2)}(\gamma)\mathcal{O}_i$.
This fact will be proved rigorously in Lemma {\ref{abelian}}.

\end{Rem}

\medskip

Now we give a precise definition of multiplication and involution on the
double cosets $K_n(\infty)\diagdown G\times G\diagup K_n(\infty)$. Denote
by $\omega_m^{(n)}$ the permutation from $\mathfrak{S}_n(\infty)$, which
acts as follows:
\begin{eqnarray*}
\omega_m^{(n)}(i)=\left\{\begin{array}{rl}
i&\textit{ if } i\leq n\textit{ or } i>n+2m
\\ i+m&\textit{ if } n<i\leq n+m
\\ i-m&\textit{ if } n+m<i\leq n+2m.
\end{array}\right.
\end{eqnarray*}

\begin{Prop}\label{pr}
Let $\tilde{g},\tilde{g}^\prime$ be double cosets and $g,g^\prime$ any
elements from $\tilde{g},\tilde{g}^\prime$, respectively. There exist
$M\left( g,g^\prime\right)\in\mathbb{N}$ such that for
$m>M\left(g,g^\prime\right)$,
$\theta_n\left(g\omega_m^{(n)}g^\prime\right)$ does not depend on the
choice of $g\in\theta_n(g)=\tilde{g}$ and
$g^\prime\in\theta_n\left(g^\prime\right)=\tilde{g}^\prime$. The
multiplication $\tilde{g}\circ\tilde{g}^\prime=
\theta_n\left(g\omega_m^{(n)}g^\prime\right)$ and the involution
$\tilde{g}^*=\theta_n\left(g^{-1}\right)$ determine a structure of
$*$-semigroup on $K_n(\infty)\diagdown G\times G\diagup K_n(\infty)$ so
that $\pi_\phi^{(2,n)}$ (see (\ref{reprcoset})) is a $*-$homomorphism.
\end{Prop}
Before proving the theorem, we discuss several auxiliary assertions.

The following statement is a generalization of Theorem 2.5 from \cite{O1}.
\begin{Lm}\label{five}
Let $\mathfrak{G}\subset\mathfrak{H}$ be discrete groups with the property:
there exists $\omega\in\mathfrak{G}$ such that for any set
$\left\{\mathfrak{g}_1,\mathfrak{g}_2,\ldots,\mathfrak{g}_p \right\}\subset
\mathfrak{G}$ $(p\in\mathbb{N})$, one can choose an element
$\mathfrak{g}\in\mathfrak{G}$ for which
$\mathfrak{H}\left\{\mathfrak{g}\mathfrak{g}_1,\mathfrak{g}\mathfrak{g}_2,
\ldots,\mathfrak{g}\mathfrak{g}_p\right\}\mathfrak{H}=\mathfrak{H}\omega
\mathfrak{H}$. If $T$ is a unitary representation of $\mathfrak{G}$ in a
Hilbert space $\mathcal{H}_T$ and $P_\mathfrak{G}$ the orthogonal
projection onto the subspace
$\mathcal{H}_T^\mathfrak{G}=\left\{\xi\in\mathcal{H}_T\big|\;
T(\mathfrak{g})\xi=\xi\textit{ for all
}\mathfrak{g}\in\mathfrak{G}\right\}$, then
$P_\mathfrak{G}=P_\mathfrak{H}T(\omega)P_\mathfrak{H}$.
\end{Lm}
\begin{proof}
Let $\xi\in\mathcal{H}_T$ and $\eta=P_\mathfrak{H}\xi$. Denote by
$\mathfrak{C}_\eta$ the closure of the set of vectors of the form
\begin{eqnarray*}
\sum\limits_{j=1}^m\alpha_j T\left(\mathfrak{g}_j\right)\eta,\textit{ where
}\mathfrak{g}_j\in\mathfrak{G},\;\alpha_1,\alpha_2,\ldots,\alpha_m\ge
0\textit{ and }\sum\limits_{j=1}^m\alpha_j=1.
\end{eqnarray*}
By our construction, $T(\mathfrak{g})\mathfrak{C}_\eta=\mathfrak{C}_\eta$.
Since there exists a unique vector $\zeta\in\mathfrak{C}_\eta$ with the
property $\|\zeta\|=
\min\left\{\|\vartheta\|\big|\vartheta\in\mathfrak{C}_\eta\right\}$, one
has $T(\mathfrak{g})\zeta=\zeta$ for all $\mathfrak{g}\in\mathfrak{G}$.
Therefore, for any $\epsilon>0$ there exist
$\alpha_1,\alpha_2,\ldots,\alpha_m\ge 0$ with the properties
\begin{eqnarray*}
\sum\limits_{j=1}^m\alpha_j=1 \;\;\textit{ and
}\;\;\left\|\sum\limits_{j=1}^m\alpha_jT\left(\mathfrak{g}_j\right)\eta-
\zeta\right\|<\frac{\epsilon}{2}.
\end{eqnarray*}
Hence, $\left\|P_\mathfrak{G}\eta-\zeta\right\|<\frac{\epsilon}{2}$ and
\begin{eqnarray}\label{seven}
\left\|\sum\limits_{j=1}^m\alpha_jT\left(\mathfrak{g}\mathfrak{g}_j\right)
\eta-P_\mathfrak{G}\eta\right\|<\epsilon\;\;\textit{ for all
}\;\;\mathfrak{g}\in\mathfrak{G}.
\end{eqnarray}
If $\mathfrak{g}$ satisfies the assumptions of the Lemma, then, using
({\ref{seven}}), we have
\begin{eqnarray*}
\left\|P_{\mathfrak{H}}T(\omega)P_{\mathfrak{H}}\eta-
P_{\mathfrak{G}}\eta\right\|<\epsilon.
\end{eqnarray*}
Since $\xi$, $\epsilon$ are chosen arbitrarily and
$\eta=P_{\mathfrak{H}}\xi$, then
$P_{\mathfrak{H}}T(\omega)P_{\mathfrak{H}}=P_{\mathfrak{G}}$.
\end{proof}
\begin{Lm}\label{wk}
For any $m,n\in\mathbb{N}$, $\mathfrak{G}=K_n(\infty)$,
$\mathfrak{H}=K_{n+m}(\infty)$,
$\omega=\left(\omega_m^{(n)},\omega_m^{(n)}\right)$ satisfy the assumptions
of Lemma {\ref{five}}. Furthermore,
$\pi_\phi^{(2)}\left(\left(\omega_m^{(n)},\omega_m^{(n)}\right)\right)$
converges to $P_n$ {\it weakly} as $m\to\infty$.
\end{Lm}
\begin{proof}
Define a permutation $\omega_{l,m}^{{n}}$ as follows:
\begin{eqnarray*}
\omega_{l,m}^{(n)}(i)=\left\{
\begin{array}{rl}
i&\textit{ if } i\le n\textit{ or }i>n+l+m
\\ i+m&\textit{ if }n<i\le n+l
\\ i-l&\textit{ if }n+l<i\le n+l+m.
\end{array}\right.
\end{eqnarray*}
A simple verification demonstrates that for $M>m$
\begin{eqnarray}\label{omg}
\omega_M^{(n)}=\omega_{m,M-m}^{(n+m)}\omega_m^{(n)}\omega_{M-m,M}^{(n+m)}.
\end{eqnarray}
Now consider an arbitrary finite collection
$\mathfrak{C}=\left\{\mathfrak{k}_1,\mathfrak{k}_2,\ldots,\mathfrak{k}_p
\right\}$ of elements from $K_n(\infty)$. There exists $M\in\mathbb{N}$
with the property $\mathfrak{C}\subset K_{n+M}\cap K_n(\infty)$. Thus
\begin{eqnarray*}
\left(\omega_M^{(n)},\omega_M^{(n)}\right)\mathfrak{C}
 \subset K_{n+M}\left(\omega_M^{(n)},\omega_M^{(n)}\right).
\end{eqnarray*}
 Hence, using ({\ref{omg}}), we obtain
 \begin{eqnarray*}
\left(\omega_m^{(n)},\omega_m^{(n)}\right) \mathfrak{g}\mathfrak{C}\subset
 K_{n+m}\left(\omega_m^{(n)},\omega_m^{(n)}\right)K_{n+m},
\textit{ where } \mathfrak{g}\in K_{n+m}.
\end{eqnarray*}
Thus the first statement is proved. The last statement follows from Lemmas
{\ref{tm}} and {\ref{five}}.
\end{proof}
\begin{proof}[{\bf Proof of Proposition \ref{pr}}]
There exists $m\in\mathbb{N}$ such that $g,g^\prime\in G_{n+m}\times
G_{n+m}$. Let $h$ be an element of $K_{n+M}\cap K_n(\infty)$, where $M\ge
n$. Notice, that $\omega_M^{(n)}h\omega_M^{(n)}\in K_{n+M}(\infty)$. Hence,
using (\ref{omg}), we have
\begin{eqnarray}
\theta_n\left(g\omega_m^{(n)}g^\prime\right)=
\theta_n\left(g\omega_M^{(n)}g^\prime\omega_M^{(n)}h\omega_M^{(n)}\right)=
\theta_n\left(gh\omega_M^{(n)}g^\prime\right)=
\theta_n\left(gh\omega_m^{(n)}g^\prime\right)
\end{eqnarray}
In a similar way, $\theta_n\left(g\omega_m^{(n)}g^\prime\right)=
\theta_n\left(g\omega_m^{(n)}hg^\prime\right)$. Thus the first statement is
proved.

Since by Lemmas {\ref{five}} and {\ref{wk}}, $g,g^\prime\in G_{n+m}\times
G_{n+m}$ and $P_n=P_{n+m}\pi_\phi^{(2)}\left(\omega_m^{(n)}\right)P_{n+m}$,
one has
\begin{eqnarray*}
P_n\pi_\phi^{(2)}\left(g\right)P_n\pi_\phi^{(2)}\left(g^\prime\right)P_n=
P_n\pi_\phi^{(2)}\left(g\right)P_{n+m}\pi_\phi^{(2)}
\left(\omega_m^{(n)}\right)P_{n+m}\pi_\phi^{(2)}\left(g^\prime\right)P_n
\\ =P_n P_{n+m}\pi_\phi^{(2)}\left(g\right)\pi_\phi^{(2)}
\left(\omega_m^{(n)}\right)\pi_\phi^{(2)}\left(g^\prime\right)P_{n+m}P_n=P_n
\pi_\phi^{(2)}\left(g\omega_m^{(n)}g^\prime\right)P_n.
\end{eqnarray*}
Finally, it is clear that
$P_n\pi_\phi^{(2)}\left(g\right)P_n^*=P_n\pi_\phi^{(2)}
\left(g^{-1}\right)P_n$.
\end{proof}

\section{A proof of the main result.}\label{mainresult} The proof of
Theorem {\ref{mainth}} splits into a few lemmas.

Define for $k\in\mathbb{N}$ the element
$\gamma(\{k\})=\left(\gamma_1(\{k\}),\gamma_2(\{k\}),\cdots,
\gamma_l(\{k\}),\cdots\right)\in\Gamma_e^\infty$ as follows:
\begin{eqnarray}\label{delta}
\gamma_l(\{k\})= \left\{
\begin{array}{rl}
\gamma_k&\text{ if }l=k
\\ e&\textit{otherwise}.
\end{array}\right.
\end{eqnarray}
For each {\it indecomposable} character $\phi$ let
$\left(\pi_\phi,\mathcal{H}_\phi,\xi_\phi\right)$ denote the cyclic
representation of the group $\mathfrak{S}_\infty\wr\Gamma$ associated to
$\phi$ via the GNS-construction.
\begin{Lm}\label{abelian}
If a $W^*-$algebra $\mathfrak{A}$ is generated by the operators
$\pi_\phi\left(\Gamma_e^\infty\right)$,
$\left\{\mathcal{O}_j\right\}_{j\in\mathbb{N}}$, and
$\mathcal{C}\left(\mathfrak{A}\right)$ is a center of $\mathfrak{A}$, then
$\left\{\mathcal{O}_j\right\}_{j\in\mathbb{N}}\subset\mathcal{C}
\left(\mathfrak{A}\right)$.
\end{Lm}
\begin{proof}
The relation
$\mathcal{O}_k\cdot\mathcal{O}_l=\mathcal{O}_l\cdot\mathcal{O}_k$ allows an
easy verification by definition ({\ref{Cesaro}}) (see {\cite{Ok1}} or
{\cite{Ok2}}).

Now prove the relation
\begin{eqnarray}\label{abel}
\mathcal{O}_l\cdot\pi_\phi\left(\gamma\right)=
\pi_\phi\left(\gamma\right)\cdot\mathcal{O}_l\textit{ for all
}\gamma\in\Gamma_e^\infty\textit{ and }l\in\mathbb{N}.
\end{eqnarray}
Let $K_n^{\mathfrak{S}}(\infty)=
K_n(\infty)\cap\left(\mathfrak{S}_\infty\times\mathfrak{S}_\infty\right)$
and $K_n^{\mathfrak{S}}(m)=K_n^{\mathfrak{S}}(\infty)\cap\left(G_m\times
G_m\right)$, where $m>n$. If $P_n^{\mathfrak{S}}$ stands for the orthogonal
projection onto $\mathcal{H}_\phi^{K_n^{\mathfrak{S}}(\infty)}$, then
\begin{eqnarray}\label{average}
P_n^{\mathfrak{S}}=\lim\limits_{m\to\infty}
\frac{1}{(m-n)!}\sum\limits_{g\in K_n^{\mathfrak{S}}(m)}
\pi_\phi^{(2)}\left(g\right)
\end{eqnarray}
in the strong operator topology and $P_n^{\mathfrak{S}}\ge P_n$. Hence,
using (\ref{Cesaro}) and (\ref{average}), we obtain for $i\le n <k$
\begin{eqnarray}
P_n^{\mathfrak{S}}\mathcal{O}_iP_n^{\mathfrak{S}}=
P_n^{\mathfrak{S}}\pi_\phi\left((i,k)\right)P_n^{\mathfrak{S}}\textit{ and
}P_n\mathcal{O}_iP_n=P_n\pi_\phi\left((i,k)\right)P_n.
\end{eqnarray}
In the case when $\gamma_l=e$ the equality (\ref{abel}) easily follows from
(\ref{Cesaro}). Therefore, it suffices to prove (\ref{abel}) for the
elements $\gamma =\gamma(\{l\})$ (see ({\ref{delta}})).

If $i\leq n <k$, then, using ({\ref{Cesaro}}), we have
\begin{eqnarray*}
& &P_n\pi_\phi\left(\gamma(\{i\})\right)\mathcal{O}_iP_n
\stackrel{\left\{P_n^{\mathfrak{S}}\ge P_n\right\}}{=}P_nP_n^{\mathfrak{S}}
\pi_\phi\left(\gamma(\{i\})\right)\mathcal{O}_iP_n^{\mathfrak{S}} P_n
\\ &\stackrel{\{({\ref{Cesaro})},{(\ref{average})}\}}{=}&P_n\pi_\phi
\left(\gamma(\{i\})\right)P_n^{\mathfrak{S}}\pi_\phi\left((i,k)\right)
P_n^{\mathfrak{S}}P_n
\\ &=&P_nP_n^{\mathfrak{S}}\pi_\phi\left((i,k)\right)\pi_\phi
\left(\gamma(\{k\})\right)P_n^{\mathfrak{S}}P_n
\\ &=&P_nP_n^{\mathfrak{S}}\pi_\phi\left((i,k)\right)\pi_\phi
\left(\gamma(\{k\})\right)\pi_\phi^{(2)}
\left(\left(\gamma(\{k\})^{-1},\gamma(\{k\})^{-1}\right)\right)P_n
\\ &\stackrel{({\ref{piphi2}})}{=}&P_nP_n^{\mathfrak{S}}\pi_\phi^{(2)}
\left(\left(e,\gamma(\{k\})^{-1}\right)\right)\pi_\phi\left((i,k)\right)P_n
\\ &=&P_n\pi_\phi^{(2)}\left(\left(\gamma(\{k\}),\gamma(\{k\})\right)\right)
\pi_\phi^{(2)}\left(\left(e,\gamma(\{k\})^{-1}\right)\right)\pi_\phi
\left((i,k)\right)P_n
\\ &=&P_n\pi_\phi\left(\gamma(\{k\})\right)\pi_\phi\left((i,k)\right)P_n=P_n
\pi_\phi\left((i,k)\right)\pi_\phi\left(\gamma(\{i\})\right)P_n
\\ &\stackrel{({\ref{Cesaro}})}{=}&P_n\mathcal{O}_i\pi_\phi
\left(\gamma(\{i\})\right)P_n.
\end{eqnarray*}
Since $\lim\limits_{n\to\infty}P_n=\mathcal{I}_{\mathcal{H}_\phi}$ (see
Lemma {\ref{tm}}), the relation
\begin{eqnarray*}
\pi_\phi\left(\gamma(\{i\})\right)\mathcal{O}_i=\mathcal{O}_i\pi_\phi
\left(\gamma(\{i\})\right)\text{ follows}.
\end{eqnarray*}
\end{proof}

We use the notation $\left(i_0,i_1,\ldots,i_{q-1}\right)$ for the cyclic
permutation $s$ which acts as follows
\begin{eqnarray*}
s\left(i\right)=\left\{
\begin{array}{rl}
i_{k+1\left(mod\,q\right)}&\text{ if }i=i_k\in
\left\{i_0,i_1,\ldots,i_{q-1} \right\}
\\ i&\textit{otherwise},
\end{array}\right.
\end{eqnarray*}

\begin{Lm}\label{power}
If $\mathcal{O}_i$ is defined as in ({\ref{Cesaro}}) and
\begin{eqnarray*}
\mathbb{D}(m,n,q)=\left\{\overrightarrow{k}=
\left(k_1,k_2,\cdots,k_q\right)\in\mathbb{N}\big|k_i\ne k_j\textit{ and }
m<k_i\leq n\,\forall i,j=1,\ldots,q\right\},
\end{eqnarray*}
then for every positive integer $m$
\begin{eqnarray*}
\mathcal{O}_i^q=\lim\limits_{n\to\infty}
\frac{1}{n^q}\sum\limits_{\overrightarrow{k}\in\mathbb{D}(m,n,q)}\pi_\phi
\left(\left(k_q,k_{q-1},\ldots,k_1,i\right) \right).
\end{eqnarray*}
\end{Lm}
\begin{proof}
If we notice that
\begin{eqnarray*}
\left(i,k_1 \right)\cdot\left(i,k_2\right)\cdots\left(i,k_q\right)=
\left(k_q,k_{q-1},\ldots,k_1,i\right)
\end{eqnarray*}
for pairwise different $i,k_1,k_2,\ldots,k_q$ and
$Card\left(\mathbb{D}(m,n)\right)=\prod\limits_{j=0}^{q-1}(n-m-j)$, the
proof becomes obvious.
\end{proof}
\begin{Lm}\label{mainlemma}
Let $g=\prod\limits_{p\in\mathbb{N}\diagup s}s_p\cdot\gamma(p)$ be a
decomposition of $g=s\cdot\gamma\in\Gamma\wr\mathfrak{S}_\infty$ (see
({\ref{decompositiontocycles}})) and $i(p)$ any element from
$p\in\mathbb{N}\diagup s$. Define $\gamma^{(i(p))}\in\Gamma_e^\infty$ as
follows
\begin{eqnarray}
\gamma_k^{(i(p))}=\left\{
\begin{array}{rl}
\gamma_{i(p)}\cdot\gamma_{s^{-1}(i(p))}\cdots\gamma_{s^{(-|p|+2)}(i(p))}
\cdot\gamma_{s^{(-|p|+1)}(i(p))}&\textit{ if }k=i(p),
\\ e&\textit{otherwise}
\end{array}\right.
\end{eqnarray}
If $\phi$ is an indecomposable character on $\Gamma\wr\mathfrak{S}_\infty$,
then
\begin{eqnarray}\label{mainrelation}
\left(\pi_\phi\left(s\cdot\gamma\right)
\prod\limits_{j}\mathcal{O}_j^{r_j}\xi_\phi,\xi_\phi\right)=
\prod\limits_{p\in\mathbb{N}\diagup
s}\left(\pi_\phi\left(\gamma^{(i(p))}\right)
\mathcal{O}_{i(p)}^{|p|-1+\sum\limits_{j\in p}r_j}\xi_\phi,\xi_\phi\right).
\end{eqnarray}
\end{Lm}
\begin{proof}
By Proposition {\ref{multiplicativity}} we have
\begin{eqnarray}
\left(\pi_\phi\left(s\cdot\gamma\right)
\prod\limits_{j}\mathcal{O}_j^{r_j}\xi_\phi,\xi_\phi\right)=
\prod\limits_{p\in\mathbb{N}\diagup
s}\left(\pi_\phi\left(s_p\cdot\gamma(p)\right)\prod\limits_{j\in
p}\mathcal{O}_j^{r_j}\xi_\phi,\xi_\phi\right).
\end{eqnarray}
Therefore it suffices to prove (\ref{mainrelation}) in the case when $s$ is
a single cycle and $\gamma=\gamma(p)$, where $p\in\mathbb{N}\diagup s$ and
$|p|>1$. Let $s=\left( i_1,i_2,\ldots,i_{|p|}\right)$. By a virtue of Lemma
{\ref{contigu}}, we find $\tilde{\gamma}\in\Gamma_e^\infty$ such that
\begin{eqnarray}
\tilde{\gamma}\cdot
s\cdot\gamma\cdot\tilde{\gamma}^{-1}=s\cdot\gamma^{(i_1)}.
\end{eqnarray}
Thus, by Lemma {\ref{abelian}},
\begin{eqnarray}
\left(\pi_\phi\left(s\cdot\gamma\right)\prod\limits_{j\in
p}\mathcal{O}_j^{r_j}\xi_\phi,\xi_\phi\right)=
\left(\pi_\phi\left(\gamma^{(i_1)}\right)\pi_\phi\left(s\right)
\prod\limits_{j\in p}\mathcal{O}_j^{r_j}\xi_\phi,\xi_\phi\right).
\end{eqnarray}
Let
\begin{eqnarray*}
\mathfrak{S}_\infty^j=
\left\{\tau\in\mathfrak{S}_\infty\big|\tau(j)=j\right\}.
\end{eqnarray*}
Now use Lemma {\ref{power}} to obtain
\begin{eqnarray*}
\begin{split}
&\left(\pi_\phi\left(\gamma^{(i_1)}\right)\pi_\phi\left(s\right)
\prod\limits_{j\in p}\mathcal{O}_j^{r_j}\xi_\phi,\xi_\phi\right)
\\ &=\lim\limits_{n\to\infty}\frac{1}{n^q}
\sum\limits_{\overrightarrow{k}\in
\mathbb{D}(m,n,q)}\Bigg(\pi_\phi\left(\gamma^{(i_1)}\right)\pi_\phi
\bigg(\Big(k_{r_{i_1}}^{\left(i_1\right)},k_{r_{i_1}-1}^{\left(i_1\right)},
\ldots,k_1^{\left(i_1\right)},i_2,
\\ & k_{r_{i_2}}^{\left(i_2\right)},\ldots,k_1^{\left(i_2\right)},i_3,
\ldots, i_{|p|},k_{r_{i_{|p|}}}^{\left(i_{|p|}\right)},\ldots,
k_1^{\left(i_{|p|}\right)},i_1\Big)\bigg)\xi_\phi,\xi_\phi\Bigg),
\end{split}
\end{eqnarray*}
where
\begin{eqnarray*}
\overrightarrow{k}=\Big(k_{r_{i_1}}^{\left(i_1\right)},
k_{r_{i_1}-1}^{\left(i_1\right)},\ldots,k_1^{\left(i_1\right)},
k_{r_{i_2}}^{\left(i_2\right)},\ldots,k_1^{\left(i_2\right)},\ldots,
k_{r_{i_{|p|}}}^{\left(i_{|p|}\right)},\ldots,
k_1^{\left(i_{|p|}\right)}\Big),\;q=\sum\limits_{j\in p}r_j.
\end{eqnarray*}
Hence, by the relation $\tau\cdot\gamma^{(i_1)}\tau^{-1}=\gamma^{(i_1)}\;
\left(\tau\in\mathfrak{S}_\infty^{i_1}\right)$, we have
\begin{eqnarray*}
\begin{split}
&\left( \pi_\phi\left(\gamma^{(i_1)}\right)\pi_\phi\left(s\right)
\prod\limits_{j\in p}\mathcal{O}_j^{r_j}\xi_\phi,\xi_\phi\right)
\\ &=\lim\limits_{n\to\infty}\frac{1}{n^{q^\prime}}
\sum\limits_{\overrightarrow{k}\in\mathbb{D}(m,n,q^\prime)}
\Bigg(\pi_\phi\left(\gamma^{(i_1)}\right)\pi_\phi
\bigg(\Big(k_{r_{i_1}}^{\left(i_1\right)},k_{r_{i_1}-1}^{\left(i_1\right)},
\ldots,k_1^{\left(i_1\right)},i_2,
\\ & k_{r_{i_2}}^{\left(i_2\right)},\ldots,k_1^{\left(i_2\right)},i_3,
\ldots,i_{|p|},k_{r_{i_{|p|}}}^{\left(i_{|p|}\right)},\ldots,
k_1^{\left(i_{|p|}\right)},i_1\Big)\bigg)\xi_\phi,\xi_\phi\Bigg),
\end{split}
\end{eqnarray*}
where
\begin{eqnarray*}
&\overrightarrow{k}=\Big(k_{r_{i_1}}^{\left(i_1\right)},
k_{r_{i_1}-1}^{\left(i_1\right)},\ldots,k_1^{\left(i_1\right)},i_2,
k_{r_{i_2}}^{\left(i_2\right)},\ldots,k_1^{\left(i_2\right)},i_3,\ldots,
i_{|p|},k_{r_{i_{|p|}}}^{\left(i_{|p|}\right)},\ldots,
k_1^{\left(i_{|p|}\right)}\Big),
\\ &q^\prime=|p|-1+\sum\limits_{j\in p}r_j.
\end{eqnarray*}
This relation, in view of Lemma {\ref{power}}, implies the statement of
Lemma {\ref{mainlemma}}.
\end{proof}
We use the notation $\mathfrak{A}_j$ for the \label{apage} $W^*-$algebra
generated by $\pi_\phi\left(\gamma\right)$,
$\gamma=\left(e,\cdots,e,\gamma_j,e,\cdots\right)$, and $\mathcal{O}_j$.
Given an operator $A$ from $\mathfrak{A}_j$, denote by $A^{(k)}$ its copy
in $\mathfrak{A}_k$:
\begin{eqnarray*}
A^{(k)}=\pi_\phi\left((j,k)\right)A\pi_\phi\left((j,k)\right)\;\;
\left(A^{(j)}=A
\right).
\end{eqnarray*}
The next assertion follows from Lemma {\ref{mainlemma}}.
\begin{Lm}\label{mainlemma1}
Let $s$, $i(p)$ be the same as in Lemma {\ref{mainlemma}}. If
$A_j,B_j\in\mathfrak{A}_j$, then
\begin{eqnarray}
\begin{split}
&\left(\pi_\phi\left(s\right)\prod\limits_jA_j\xi_\phi,
\prod\limits_jB_j\xi_\phi\right)
\\ =&\prod\limits_{p\in \mathbb{N}\diagup s}\bigg(A_{i(p)}^{(i(p))}
\left(B_{i(p)}^{(i(p))}\right)^*A_{s^{-1}(i(p))}^{(i(p))}
\left(B_{s^{-1}(i(p))}^{(i(p))}\right)^*\cdots
\\ &\cdots A_{s^{1-|p|}(i(p))}^{(i(p))}
\left(B_{s^{1-|p|}(i(p))}^{(i(p))}\right)^*
\mathcal{O}_{i(p)}^{|p|-1}\xi_\phi,\;\xi_\phi\bigg)
\end{split}
\end{eqnarray}
\end{Lm}
The following lemma is an analogue of Theorem 1 from {\cite{Ok2}}.
\begin{Lm}\label{discrete}
Let $\Delta=[a,b]$ be an interval in $[-1,0]$ or in $[0,1]$ with the
property $\min\left\{|a|,|b|\right\}>\varepsilon>0$. If $E_{\Delta}^{(i)}$
is a spectral projection of $\mathcal{O}_i$ corresponding to $\Delta$, then
for any orthogonal projection $E$ from $\mathfrak{A}_i$ one has
$\left(EE_{\Delta}^{(i)}\xi_\phi,\xi_\phi\right)^2\ge\varepsilon
\left(EE_{\Delta}^{(i)}\xi_\phi,\xi_\phi\right)$.
\end{Lm}
\begin{proof}
Using Lemmas {\ref{abelian}} and {\ref{mainlemma1}}, we have
\begin{eqnarray}\label{est}
\begin{split}
&\left|\left(\pi_\phi\left((i,i+1)\right)EE_{\Delta}^{(i)}\xi_\phi,
EE_{\Delta}^{(i)}\xi_\phi\right)\right|&
\\ &=\left|\left(\mathcal{O}_iEE_{\Delta}^{(i)}\xi_\phi,
EE_{\Delta}^{(i)}\xi_\phi\right)\right|>\varepsilon\left|
\left(EE_{\Delta}^{(i)}\xi_\phi,\xi_\phi\right)\right|.&
\end{split}
\end{eqnarray}
On the other hand, under the assumption
$E^{(i+1)}=\pi_\phi\left((i,i+1)\right)E\pi_\phi\left((i,i+1)\right)$, one
has
\begin{eqnarray*}
EE_{\Delta}^{(i)}\cdot E^{(i+1)}E_{\Delta}^{(i+1)}\cdot\pi_\phi
\left((i,i+1)\right)
\\ =\pi_\phi\left((i,i+1)\right)\cdot EE_{\Delta}^{(i)}\cdot E^{(i+1)}
E_{\Delta}^{(i+1)}.
\end{eqnarray*}
Therefore,
\begin{eqnarray*}
&\left|\left(\pi_\phi\left((i,i+1)\right)EE_{\Delta}^{(i)}\xi_\phi,
EE_{\Delta}^{(i)}\xi_\phi\right)\right|
\\ &=\left|\left(\pi_\phi\left((i,i+1)\right)E^{(i+1)}E_{\Delta}^{(i+1)}
EE_{\Delta}^{(i)}\xi_\phi,\xi_\phi\right)\right|
\\ &\le\left|\left(E^{(i+1)}E_{\Delta}^{(i+1)}EE_{\Delta}^{(i)}\xi_\phi,
\xi_\phi\right)\right|\stackrel{(Prop. \ref{multiplicativity})}{=}
\left(EE_{\Delta}^{(i)}\xi_\phi,\xi_\phi\right)^2.
\end{eqnarray*}
Hence, using ({\ref{est}}), we obtain our statement.
\end{proof}
The following statement is well known (see {\cite{Ok2}}) and also follows
from Lemma {\ref{discrete}}.
\begin{Co}
There exists at most countable set of numbers $\alpha_i$, $\beta_i$ from
$(0,1)$ and a set of pairwise orthogonal projections
$\left\{E^{(k)}\left(\alpha_i\right),E^{(k)}\left(\beta_i\right)\right\}
\subset\mathfrak{A}_k$ such that
\begin{eqnarray}\label{spectraldecomposition}
\mathcal{O}_k=\sum\alpha_iE^{(k)}\left(\alpha_i\right)-\sum\beta_iE^{(k)}
\left(\beta_i\right).
\end{eqnarray}
\end{Co}
The following assertion is an analogue of Theorem 2 from {\cite{Ok2}}.
\begin{Lm}\label{integer}
Let $r$ be a number from $\left\{\alpha_i,\beta_i\right\}$ and let $E$ be
any projection from $\mathfrak{A}_k$. If $\left(E\cdot
E^{(k)}(r)\xi_\phi,\xi_\phi\right)=r\nu(r)\ne 0$, then
$\nu(r)\in\mathbb{Z}$.
\end{Lm}
\begin{proof}
For completeness of the proof, we use the arguments of Kerov, Olshanski,
Vershik and Okounkov from {\cite{KOV}} and {\cite{Ok2}}.

For any $m\in\mathbb{N}$, define the projection $e_m(r)$ as follows:
\begin{eqnarray*}
&e_m(r)=\prod\limits_{j=1}^m E^{(j)}\cdot E^{(j)}(r),\textit{ where }
\\ & E^{(j)}=\pi_\phi\left((j,k)\right)E\pi_\phi\left((j,k)\right),\;
E^{(j)}(r)=\pi_\phi\left((j,k)\right)E^{(k)}(r)\pi_\phi\left((j,k)\right).
\end{eqnarray*}
Let $\mathbb{P}_m(s)$ be the set of orbits $s$ on $\{1,2,\ldots,m\}$. If
$s\in\mathfrak{S}_m$, then by Lemma {\ref{mainlemma1}} we obtain
\begin{eqnarray}\label{orb}
\left(\pi_\phi(s)e_m(r)\xi_\phi,\;e_m(r)\xi_\phi\right)=
\nu(r)^{\left|\mathbb{P}_m(s)\right|}\prod\limits_{p\in\mathbb{P}_m(s)}
r^{|p|}.
\end{eqnarray}
Set
$\phi_r(s)=\dfrac{\left(\pi_\phi(s)e_m(r)\xi_\phi,\;e_m(r)\xi_\phi\right)}
{\left(e_m(r)\xi_\phi,\;e_m(r)\xi_\phi\right)}$. Using ({\ref{orb}}), we
have
\begin{equation}\label{phir}
\phi_r(s)=\dfrac{\nu(r)^{\left|\mathbb{P}_m(s)\right|}}{\nu(r)^m}.
\end{equation}
Therefore, $\phi_r$ is an indecomposable character on $\mathfrak{S}_\infty$
in view of Proposition {\ref{multiplicativity}}.

We following G. Olshanski (see {\cite{O1}}) in expounding the proof of the
following formula:
\begin{equation}\label{alt}
\sum\limits_{s\in\mathfrak{S}_m}sgn(s)\;t^{\left|\mathbb{P}_m(s)\right|}=
t(t-1)\cdots(t-m+1).
\end{equation}
For that, we consider the canonical projection $p_{m,m-1}$ from
$\mathfrak{S}_m$ onto $\mathfrak{S}_{m-1}$:
\begin{eqnarray*}
\left(p_{m,m-1}(s)\right)(i)=\left\{
\begin{array}{rl}
s(i)&\textit{ if } s(i)<m
\\ s(m)&\textit{ if } s(i)=m.
\end{array}\right.
\end{eqnarray*}
Since $\left|\mathbb{P}_{m-1}\left(p_{m,m-1}(s)\right)\right|=
\left|\mathbb{P}_m(s)\right|$ when $s\notin \mathfrak{S}_{m-1}$, and
$\left|\mathbb{P}_{m-1}\left(p_{m,m-1}(s)\right)\right|=
\left|\mathbb{P}_m(s)\right|-1$ when $s\in\mathfrak{S}_{m-1}$, then
\begin{eqnarray*}
\sum\limits_{s\in\mathfrak{S}_m}sgn(s)\;t^{\left|\mathbb{P}_m(s)\right|}=
\sum\limits_{s\in\mathfrak{S}_{m-1}}
\sum\limits_{\tilde{s}\in\mathfrak{S}_m:\;p_{m,m-1}(\tilde{s})=s}sgn(s)\;
t^{\left|\mathbb{P}_m(s)\right|}=
\\ t\cdot\sum\limits_{s\in\mathfrak{S}_{m-1}}
t^{\left|\mathbb{P}_m(s)\right|}-(m-1)\cdot
\sum\limits_{s\in\mathfrak{S}_{m-1}}t^{\left|\mathbb{P}_m(s)\right|}=(t-m+1)
\sum\limits_{s\in\mathfrak{S}_{m-1}}t^{\left|\mathbb{P}_m(s)\right|}.
\end{eqnarray*}
Hence {(\ref{alt})} is now accessible by an elementary induction argument.

We follow the idea of A. Okounkov in considering the orthogonal projection
\begin{eqnarray*}
Alt_r(m)=\frac{1}{m!}\sum\limits_{s\in\mathfrak{S}_m}sgn(s)\;
\pi_{\phi_r}(s).
\end{eqnarray*}
Since $\sum\limits_{s\in\mathfrak{S}_m}sgn(s)\;{\phi_r}(s)\ge 0$, then,
using ({\ref{phir}}) and ({\ref{alt}}), we obtain for $r>0$
\begin{eqnarray*}
\nu(r)\cdot\left(\nu(r)-1\right)\cdots\left(\nu(r)-m+1\right)\ge 0\textit{
for all }m\in\mathbb{N}.
\end{eqnarray*}
Thus, we get a contradiction in the case $\nu(r)>0$. The opposite case
$\nu(r)<0$ can be considered in a similar way. For that, one should use the
formula
\begin{equation*}
\sum\limits_{s\in\mathfrak{S}_m}t^{\left|\mathbb{P}_m(s)\right|}=
t(t+1)\cdots(t+m-1)\textit{ (see \cite{O1})}
\end{equation*}
and consider the projection
\begin{eqnarray*}
Sym_r(m)=\frac{1}{m!}\sum\limits_{s\in\mathfrak{S}_m}\pi_{\phi_r}(s).
\end{eqnarray*}
\end{proof}
\begin {proof}[\bf{Proof of Theorem \ref{mainth}}]
Let $E_k(r)$ be the spectral projection of $\mathcal{O}_k$ (see
(\ref{Cesaro}),(\ref{spectraldecomposition})). By Lemma {\ref{integer}},
for $r\ne 0$ the $W^*-$algebra $E_k(r)\mathfrak{A}_k$ (see p.
\pageref{apage}) is finite-dimensional. On the other hand, use Lemma
{\ref{abelian}} to obtain the unitary representation
$\left(E_k(r)\pi_\phi\Big|_\Gamma, E_k(r)\mathcal{H}_\phi\right)$ of the
group $\Gamma$ in the space $E_k(r)\mathcal{H}_\phi$. Thus, the
representations $\varrho^r$ for $r\ne 0$ as in Theorem {\ref{mainth}} are
the irreducible components of $\left(E_k(r)\pi_\phi\Big|_\Gamma,
E_k(r)\mathcal{H}_\phi\right)$. The formula for characters follows from
Lemmas {\ref{abelian}} and {\ref{mainlemma1}}. Finally, for each character
as in Theorem {\ref{mainth}} we construct the realization as in Section
{\ref{two}}.
\end{proof}

Authors:

\medskip

Nessonov Nikolay, Institute For Low Temperature Physics and Engineering,
Department of Mathematics, 47 Lenin Avenue, Kharkiv, Ukraine,
(0572)-30-85-85

 nessonov@ilt.kharkov.ua

\medskip

Dudko Artem, Kharkov National University,

artemdudko@rambler.ru
\end{document}